\documentclass[12pt]{amsproc}
\usepackage{amsmath,amssymb,amscd}

\begin{document}

\title[Free Probability for Pairs of Faces I]{Free Probability for Pairs of Faces I}
\author{Dan-Virgil Voiculescu}
\address{D.V. Voiculescu \\ Department of Mathematics \\ University of California at Berkeley \\ Berkeley, CA\ \ 94720-3840}
\thanks{Research supported in part by NSF Grant DMS-1301727.}
\date{}

\begin{abstract}
We consider a notion of bi-freeness for systems of non-commutative random variables with two faces, one of left variables and another of right variables. This includes bi-free convolution operations, bi-free cumulants and the bi-free central limit.

\end{abstract}

\maketitle

\setcounter{section}{-1}
\section{Introduction}
\label{sec0}

On a free product of Hilbert spaces with specified unit vector there are two actions of the operators of the initial spaces, corresponding to a left and to a right tensorial factorization respectively (\cite{12}, \cite{16}). In usual free probability this gives a choice between left and right when modeling free random variables on free products of spaces. From here it is then natural to imagine the situation when algebras of left operators and algebras of right operators on the free product space are considered at the same time. We will show that there is a ``two-faced'' extension of free probability to deal with the expectation values for such combined systems of operators and begin developing in this paper and its successors the basics of the theory. Like with usual free probability we have found it convenient in the development of basic concepts to move between a purely algebraic framework and a $C^*$-algebra functional analysis framework.

We show in what follows that there is a notion of ``bi-free independence'' (or ``bi-freeness'') for two-faced systems which has the right properties for the development of probability theory concepts based on it. Restricted to systems which have only left (or only right) non-constant variables, bi-freeness reduces to freeness, however when there are both left and right variables, bi-freeness implies for instance that certain groups of variables at the level of their distributions exhibit commutative (``tensor product'') independence.  Bi-freeness may be related to a number of recent non-commutative independencies (\cite{2}, \cite{1}, \cite{5}, \cite{6}, \cite{7}, \cite{11}). An exploration of possible connections would be of much interest.

Addition and multiplication of bi-free two-faced families of non-commutative random variables give rise to additive and multiplicative bi-free convolution operations for the distributions. We prove, based on general Lie theory considerations, the existence of bi-free cumulants which linearize the additive bi-free convolution. Like in our earliest paper on free probability \cite{12}, these general nonsense results are sufficient in order to prove a central limit theorem. We are also able to identify the bi-free central limit distributions, that is the bi-free analogue of Gaussian distributions. They arise from left and right creation and annihilation operators on the full Fock space. (For other recent central limit work related to free probability see \cite{3}, \cite{4}.)

This paper has eight sections including this introduction. Section~\ref{sec1} is devoted to preliminaries on the free product of vector-spaces with specified state vector, which is an essential technical ingredient for the rest of the paper. In Section~\ref{sec2} we introduce two-faced probabilistic notions and show that bi-freeness has the right properties for a notion of independence in this framework. These considerations are purely algebraic.

Then, in the short Section~\ref{sec3} we discuss additional features when we deal with two-faced $C^*$-algebraic systems and bi-freeness.

Section~\ref{sec4} introduces bi-free convolution operations, the operations on distributions arising from operations on bi-free two-faced families of non-commutative random variables.

Bi-free cumulants are the subject of Section~\ref{sec5}, the main result being an abstract existence theorem.

A few basic examples of bi-freeness are discussed in Section~\ref{sec6}.

Using bi-free cumulants and the left and right creation and annihilation operators example we find in Section~\ref{sec7} the central limit distributions and prove an algebraic central limit theorem. We consider a general case as well as the hermitian and the $*$-distributions case.

If $B$ is an algebra over ${\mathbb C}$, we sketch a definition of bi-freeness over $B$ in Section~\ref{sec8}, together with a definition of an assorted $B'-B$ non-commutative probability space, to replace the usual non-commutative probability space over $B$ (\cite{17}).

Finally, our choice to call the sets of ``left variables'' and ``right variables'' a ``pair of faces'', was inspired by the two faces of the mythical Janus, though it is open whether there are mathematical examples, where one set of variables is looking into the future, while the other is looking into the past, like the two faces of Janus. We were also motivated to avoid calling the systems ``bipartite'', since ``bipartite'' will be reserved for the particular class of systems where the sets of left and right variables commute with each other.

\section{Preliminaries on Free Products of Vector Spaces}
\label{sec1}

\noindent
{\bf 1.1.} We collect here the basics about left and right operators on a free product of Hilbert spaces and we also give similar purely algebraic concepts, which we will use later.

\bigskip
\noindent
{\bf 1.2.} Let $({\mathcal H}_i,\xi_i)_{i \in I}$ be a family of Hilbert spaces with specified unit vectors. The Hilbert space free product $({\mathcal H},\xi) = \underset{i \in I}{*} ({\mathcal H}_i,\xi_i)$ is given by
\[
{\mathcal H} = {\mathbb C}\xi \oplus \underset{n \ge 1}{\bigoplus} \left( \underset{i_1 \ne i_2 \ne \dots \ne i_n}{\bigoplus} \overset{\circ}{\mathcal H}_{i_1} \otimes \dots \otimes \overset{\circ}{\mathcal H}_{i_n}\right)
\]
where $\overset{\circ}{\mathcal H}_i = {\mathcal H}_i \ominus {\mathbb C}\xi_i$. Here the tensor products and direct sums are hilbertian and involve completions and orthogonality and $\xi$ is given norm~$1$.

\bigskip
\noindent
{\bf 1.3.} If $i \in I$ we define
\[
{\mathcal H}(l,i) = {\mathbb C}\xi \oplus \underset{n \ge 1}{\bigoplus} \left( \underset{\underset{i_i \ne i}{i_1\ne \dots \ne i_n}}{\bigoplus} \overset{\circ}{\mathcal H}_{i_1} \otimes \dots \otimes \overset{\circ}{\mathcal H}_{i_n}\right)
\]
and
\[
{\mathcal H}(r,i) = {\mathbb C}\xi \oplus \underset{n \ge 1}{\bigoplus} \left( \underset{\underset{i_n \ne i}{i_n \ne \dots \ne i_n}}{\bigoplus} \overset{\circ}{\mathcal H}_{i_1} \otimes \dots \otimes \overset{\circ}{\mathcal H}_{i_n}\right).
\]

There are natural identifications given by unitary operators $V_i: {\mathcal H}_i \otimes {\mathcal H}(l,i) \to {\mathcal H}$ and $W_i: {\mathcal H}(r,i) \otimes {\mathcal H}_i \to {\mathcal H}$. In the case of $V_i$ the identifications are the obvious ones for
\[
\begin{aligned}
&\xi_i \otimes \xi \to \xi \\
&\overset{\circ}{\mathcal H}_i \otimes \xi \to \overset{\circ}{\mathcal H}_i \\
&\xi_i \otimes (\overset{\circ}{\mathcal H}_{i_1} \otimes \dots \otimes \overset{\circ}{\mathcal H}_{i_n}) \to \overset{\circ}{\mathcal H}_{i_1} \otimes \dots \otimes \overset{\circ}{\mathcal H}_{i_n} \\
&\overset{\circ}{\mathcal H}_i \otimes (\overset{\circ}{\mathcal H}_{i_1} \otimes \dots \otimes \overset{\circ}{\mathcal H}_{i_n}) \to \overset{\circ}{\mathcal H}_i \otimes \overset{\circ}{\mathcal H}_{i_1} \otimes \dots \otimes \overset{\circ}{\mathcal H}_{i_n}.
\end{aligned}
\]
Similarly for the $W_i$ the list of identifications is
\[
\begin{aligned}
&\xi \otimes \xi_i \to \xi \\
&\xi \otimes \overset{\circ}{\mathcal H}_i \to \overset{\circ}{\mathcal H}_i \\
&(\overset{\circ}{\mathcal H}_{i_1} \otimes \dots \otimes \overset{\circ}{\mathcal H}_{i_n}) \otimes \xi_i \to \overset{\circ}{\mathcal H}_{i_1} \otimes \dots \otimes \overset{\circ}{\mathcal H}_{i_n} \\
&(\overset{\circ}{\mathcal H}_{i_1} \otimes \dots \otimes \overset{\circ}{\mathcal H}_{i_n}) \otimes \overset{\circ}{\mathcal H}_i \to \overset{\circ}{\mathcal H}_{i_1} \otimes \dots \otimes \overset{\circ}{\mathcal H}_{i_n} \otimes \overset{\circ}{\mathcal H}_i.
\end{aligned}
\]

\bigskip
\noindent
{\bf 1.4.} If $T \in {\mathcal B}({\mathcal H}_i)$ we define the left and right operators
\[
\begin{aligned}
\lambda_i(T) &= V_i(T \otimes I_{{\mathcal H}(r,i)})V_i^{-1} \\
\rho_i(T) &= W_i(I_{{\mathcal H}(r,i)} \otimes T)W_i^{-1}.
\end{aligned}
\]
Both $\lambda_i$ and $\rho_i$ are representations of $B({\mathcal H}_i)$ on $B({\mathcal H})$.

\bigskip
\noindent
{\bf 1.5.} If $T_1 \in {\mathcal H}_i$ and $T_2 \in {\mathcal H}_{i'}$ we have the commutation relation
\[
[\lambda_i(T_1),\rho_{i'}(T_2)] = \begin{cases}
0 &\text{if $i \ne i'$} \\
[T_1,T_2] \oplus {\mathcal O}_{{\mathcal H} \ominus {\mathcal H}_i} &\text{if $i = i'$}
\end{cases}
\]
where ${\mathcal H}_i$ is identified with
\[
V_i({\mathcal H}_i \otimes \xi) = W_i(\xi \otimes {\mathcal H}_i) = {\mathbb C}\xi \oplus \overset{\circ}{\mathcal H}_i.
\]

If $i \ne i'$ this follows from an identification of $\lambda_i(T_1),\rho_{i'}(T_2)$, ${\mathcal H}$ with $T_1 \otimes I \otimes I$, $I \otimes I \otimes T_2$ and
\[
{\mathcal H}_i \otimes {\mathcal H}(lr,i,i') \otimes {\mathcal H}_{i'}
\]
where
\[
{\mathcal H}(lr,i,i') = {\mathbb C}\xi \oplus \underset{n \ge 1}{\bigoplus} \left( \underset{\underset{i_1\ne \dots \ne i_n}{i_1 \ne i,i_n \ne i'}}{\bigoplus} \overset{\circ}{\mathcal H}_{i_1} \otimes \dots \otimes \overset{\circ}{\mathcal H}_{i_n}\right).
\]

In case $i = i'$ one identifies ${\mathcal H}$ with
\[
{\mathcal H}_i \oplus ({\mathcal H}_i \otimes \overset{\circ}{\mathcal H}(lr,i) \otimes {\mathcal H}_i)
\]
where
\[
\overset{\circ}{\mathcal H}(lr,i) = \bigoplus_{n \ge 1} \left( \underset{\underset{i_1 \ne i \ne i_n}{i_1 \ne \dots \ne i_n}}{\bigoplus} \overset{\circ}{\mathcal H}_{i_1} \otimes \dots \otimes \overset{\circ}{\mathcal H}_{i_n}\right)
\]
and $\lambda_i(T_1),\rho_i(T_2)$ with $T_1 \oplus (T_1 \otimes I \otimes I)$ and $T_2 \oplus (I \otimes I \otimes T_2)$.

\bigskip
\noindent
{\bf 1.6.} There is also a purely algebraic analogue of the free product construction. For this the Hilbert space with specified unit vector $({\mathcal H},\xi)$ is replaced by a {\em vector space with specified state-vector}, which is a triple $({\mathcal X},\overset{\circ}{\mathcal X},\xi)$, where ${\mathcal X}$ is a vector space, $\overset{\circ}{\mathcal X} \subset {\mathcal X}$ is a subspace of codimension~$1$ and $0 \ne \xi \in {\mathcal X}$ is a vector so that ${\mathbb C}\xi + \overset{\circ}{\mathcal X} = {\mathcal X}$. We will often write ${\mathcal X} = {\mathbb C}\xi \oplus \overset{\circ}{\mathcal X}$, identifying ${\mathcal X}$ with the direct sum. There are {\em two other descriptions of a vector space with specified state-vector: as $({\mathcal X},p,\xi)$ and $({\mathcal X},\psi,\xi)$}. In the first, $p: {\mathcal X} \to {\mathcal X}$ is the idempotent so that $p(\xi) = \xi$ and $\ker p = \overset{\circ}{\mathcal X}$. In the second, $\psi$ is the functional $\psi: {\mathcal X} \to {\mathbb C}$ so that $\psi(\xi) = 1$ and $\ker \psi = \overset{\circ}{\mathcal X}$. Clearly, we have $p(x) = \psi(x)\xi$. We shall refer to both descriptions also as a vector space with state-vector.

Sometimes the vector $\psi$ in $({\mathcal X},p,\xi)$ will play the role of a coordinate and its choice in $p({\mathcal X})$ will only be to facilitate computations, while the essential data is the pair $({\mathcal X},p)$ where $p: {\mathcal X} \to {\mathcal X}$ is an idempotent of rank~$1$. We shall refer to $({\mathcal X},p)$ as a {\em vector space with specified vector-state}.

\bigskip
\noindent
{\bf 1.7.} Given a vector space ${\mathcal X}$, we shall denote by ${\mathcal L}({\mathcal X})$ the algebra of linear operators on ${\mathcal X}$. In the case of a vector space with specified vector-state, $({\mathcal X},p)$, there is a linear functional $\varphi: {\mathcal L}({\mathcal X}) \to {\mathbb C}$ so that $\varphi(T)p = pTp$. Since $\varphi(1) = 1$ the pair $({\mathcal L}({\mathcal X}),\varphi)$ is a non-commutative probability space. If the specification is necessary, we will use the notation $\varphi_p$.

\bigskip
\noindent
{\bf 1.8.} Given a family $({\mathcal X}_i,\overset{\circ}{\mathcal X}_i,\psi_i)$, $i \in I$, {\em the free product}
\[
({\mathcal X},\overset{\circ}{\mathcal X},\xi) = \underset{i \in I}{*} ({\mathcal X}_i,\overset{\circ}{\mathcal X}_i,\xi_i)
\]
is defined by the formulae
\[
\begin{aligned}
{\mathcal X} &= \overset{\circ}{\mathcal X} \oplus {\mathbb C}\xi, \\
\overset{\circ}{\mathcal X} &= \underset{n \ge 1}{\bigoplus} \left( \underset{i_1 \ne \dots \ne i_n}{\bigoplus} \overset{\circ}{\mathcal X}_{i_1} \otimes \dots \otimes \overset{\circ}{\mathcal X}_{i_n}\right).
\end{aligned}
\]
Note that in these formulae which look like the hilbertian ones, the $\oplus$ and $\otimes$ are purely algebraic and do not involve completions like in the case of Hilbert spaces.

Since clearly the vector space with specified vector-state $({\mathcal X},p)$ is determined by the $({\mathcal X}_i,p_i)$ we will write
\[
({\mathcal X},p) = \underset{i \in I}{*} ({\mathcal X}_i,p_i).
\]

\bigskip
\noindent
{\bf 1.9.} The same kind of formulae used in the case of Hilbert spaces also give in the purely algebraic context of a family $({\mathcal X}_i,\overset{\circ}{\mathcal X}_i,\xi_i)$, $i \in I$, definitions of ${\mathcal X}(l,i)$, ${\mathcal X}(r,i)$, $V_i$, $W_i$ and $\lambda_i: {\mathcal L}({\mathcal X}_i) \to {\mathcal L}({\mathcal X})$, $\rho_i: {\mathcal L}({\mathcal X}_i) \to {\mathcal L}({\mathcal X})$ where
\[
({\mathcal X},\overset{\circ}{\mathcal X},\xi) = \underset{i \in I}{*} ({\mathcal X}_i,\overset{\circ}{\mathcal X}_i,\xi_i).
\]

If $T \in {\mathcal L}({\mathcal X}_i)$ and $x_j \in \overset{\circ}{\mathcal X}_{i_j}$, $i_j \ne i_{j+1}$, $1 \le j \le n$, and if we put
\[
\lambda_i(T)x_1 \otimes \dots \otimes x_n = \eta,
\]
then, if $i \ne i_1$ we have
\[
\eta = (\varphi_{p_i}(T)x_1 \otimes \dots \otimes x_n) \oplus ((T\xi_i - \varphi_{p_i}(T)\xi_i) \otimes x_1 \otimes \dots \otimes x_n)
\]
and if $i = i_1$, we have
\[
\eta = ((Tx_1 - p_i(Tx_1)) \otimes x_2 \otimes \dots \otimes x_n) \oplus \psi_i(Tx_1)x_2 \otimes \dots \otimes x_n.
\]
We also have
\[
\lambda_i(T)\xi = \varphi_{p_i}(T)\xi \oplus (T\xi_i - p_i(T\xi_i))
\]
where the second component is in $\overset{\circ}{\mathcal X}_i$.

Similarly, if we put
\[
\rho_i(T)x_1 \otimes \dots \otimes x_n = \eta
\]
then, if $i \ne i_n$ we have
\[
\eta = (\varphi_{p_i}(T)x_1 \otimes \dots \otimes x_n) \oplus (x_1 \otimes \dots \otimes x_n \otimes (T\xi_i - (p_i(T\xi)))
\]
and if $i = i_n$
\[
\eta = (x_1 \otimes \dots \otimes x_{n-1} \otimes (Tx_n - p_i(Tx_n))) \oplus (\psi_i(Tx_n)x_1 \otimes \dots \otimes x_{n-1}).
\]

Also $\lambda_i(T)\xi = \rho_i(T)\xi$.

The commutation formula for $[\lambda_i(T_1),\rho_{i'}(T_2)]$ is the obvious analogue of the one in the Hilbert space case in $1.5$.

\bigskip
\noindent
{\bf 1.10.} Morphisms between vector spaces with specified state-vector $({\mathcal X},\overset{\circ}{\mathcal X},\xi)$ and $({\mathcal X}',\overset{\circ}{\mathcal X}{}',\xi')$ are linear maps $S: {\mathcal X} \to {\mathcal X}'$, so that $S\xi = \xi'$ and $S(\overset{\circ}{\mathcal X}) \subset \overset{\circ}{\mathcal X}{}'$. Clearly all relevant information about $S$ is in $\overset{\circ}{S}: \overset{\circ}{\mathcal X} \to \overset{\circ}{\mathcal X}{}'$ defined as $\overset{\circ}{S} = {}_{\overset{\circ}{\mathcal X}{}'}|S|_{\overset{\circ}{\mathcal X}}$.

Given $({\mathcal X}_i,\overset{\circ}{\mathcal X}_i,\xi_i)$, $({\mathcal X}'_i,\overset{\circ}{\mathcal X}{}'_i,\xi'_i)$ and morphisms $S_i$, $i \in I$ there is a free product morphism $S = \underset{i \in I}{*} S_i$ from $\underset{i \in I}{*} ({\mathcal X}_i,\overset{\circ}{\mathcal X}_i,\xi_i)$ to $\underset{i \in I}{*} ({\mathcal X}'_i,\overset{\circ}{\mathcal X}{}'_i,\xi'_i)$ where $S\xi = \xi'$ and
\[
S|_{\overset{\circ}{\mathcal X}_{i_1} \otimes \dots \otimes \overset{\circ}{\mathcal X}_{i_n}} = \overset{\circ}{S}_{i_1} \otimes \dots \otimes \overset{\circ}{S}_{i_n}
\]
where $i_1 \ne i_2 \ne \dots \ne i_n$.

When studying intertwining properties of $S = \underset{i \in I}{*} S_i$ it will be necessary to distinguish between the $V_i$ and $W_i$ associated with the ${\mathcal X}_i$ and ${\mathcal X}'_i$, which will be achieved by abusing notations and writing $V'_i$, $W'_i$ for the maps associated with the ${\mathcal X}'_i$. In the same vein we may add primes and write $p'$, $\lambda'_i$, $\rho'_i$, etc.

\bigskip
\noindent
{\bf 1.11. Lemma.} {\em Let $S$ be a morphism from $({\mathcal X},\overset{\circ}{\mathcal X},\xi)$ to $({\mathcal X}',\overset{\circ}{\mathcal X}',\xi')$ and $T \in {\mathcal L}({\mathcal X})$, $T' \in {\mathcal L}({\mathcal X}')$ such that $T'S = ST$. Then we have
\[
\varphi(T) = \varphi'(T').
\]
}

\bigskip
\noindent
{\bf {\em Proof.}} Indeed, since $Sp = p'S \ne 0$ we have
\[
\varphi'(T')p'S = p'T'p'S = SpTp = \varphi(T)Sp = \varphi(T)p'S
\]
which implies $\varphi(T) = \varphi'(T')$. \qed

\bigskip
\noindent
{\bf 1.12. Lemma.} {\em Let $S_i$ be morphisms from $({\mathcal X}_i,\overset{\circ}{\mathcal X}_i,\xi_i)$ to $({\mathcal X}'_i,\overset{\circ}{\mathcal X}{}'_i,\xi'_i)$ and $S = \underset{i \in I}{*} S_i$. Then we have $S\xi = \xi'$ and $S{\mathcal X}(l,i) \subset {\mathcal X}'(l,i)$, $S{\mathcal X}(r,i) \subset {\mathcal X}'(r,i)$. Moreover we have
\[
\begin{aligned}
SV_i &= V'_i(S_i \otimes S(l,i)) \text{ and} \\
SW_i &= W'_i(S(r,i) \otimes S_i)
\end{aligned}
\]
where $S(l,i) = {}_{{\mathcal X}'(l,i)}|S|_{{\mathcal X}(l,i)}$ and $S(n,i) = {}_{{\mathcal X}'(r,i)}|S|_{{\mathcal X}(r,i)}$.
}

\bigskip
The proof follows immediately from the formulae defining $V_i$, $W_i$, $V'_i$, $W'_i$ and the fact that $S_i\xi_i = \xi'_i$, $S_i\overset{\circ}{\mathcal X}_i \subset \overset{\circ}{\mathcal X}{}'_i$ and is left to the reader.

\bigskip
\noindent
{\bf 1.13. Lemma.} {\em Let $S_i$ be morphisms from $({\mathcal X}_i,\overset{\circ}{\mathcal X}_i,\xi_i)$ to $({\mathcal X}'_i,\overset{\circ}{\mathcal X}{}'_i,\xi'_i)$ and $S = \underset{i \in I}{*} S_i$. Then if $T_i \in {\mathcal L}({\mathcal X}_i)$, $T'_i \in {\mathcal L}({\mathcal X}'_i)$ are such that $T'_iS_i = S_iT_i$ we have
\[
\begin{aligned}
\lambda'_i(T'_i)S &= S\lambda_i(T_i) \text{ and} \\
\rho'_i(T'_i)S &= S\rho_i(T_i).
\end{aligned}
\]
}

\bigskip
\noindent
{\bf {\em Proof.}} This is a consequence of the preceding lemma. Indeed, we have
\[
\begin{aligned}
S\lambda_i(T_i) &= SV_i(T_i \otimes I_{{\mathcal X}(l,i)})V_i^{-1} \\
&= V'_i(S_i \otimes S(l,i))(T_i \otimes I_{{\mathcal X}(l,i)})V_i^{-1} \\
&= V'_i(T'_i \otimes I_{{\mathcal X}'(l,i)})(S_i \otimes S(l,i))V_i^{-1} \\
&= V'_i(T'_i \otimes I_{{\mathcal X}'(l,i)})V'{}^{-1}_iS \\
&= \lambda'_i(T'_i)S
\end{aligned}
\]
and a similar computation proves the intertwining of $\rho_i(T_i)$ and $\rho'_i(T'_i)$. \qed

\bigskip
\noindent
{\bf 1.13. Remark.} We record here some properties of the free product of vector spaces which the reader will find useful in the discussion of bi-freeness in \S$\ref{sec2}$. Let $(({\mathcal X}_k,\overset{\circ}{\mathcal X}_k,\xi_k))_{k \in K}$ be a family of vector spaces with specified state-vector and let $\underset{i \in I}{\coprod} K_i = K$ be a partition of the index set. Let $(({\mathcal Y}_i,\overset{\circ}{\mathcal Y}_i,\eta_i))_{i \in I}$ be the family of free products $({\mathcal Y}_i,\overset{\circ}{\mathcal Y}_i,\eta_i) = \underset{k \in K_i}{*} ({\mathcal X}_k,\overset{\circ}{\mathcal X}_k,\xi_k)$. Then there are natural identifications $\sigma: \underset{i \in I}{*} ({\mathcal Y}_i,\overset{\circ}{\mathcal Y}_i,\eta_i) \to \underset{k \in K}{*} ({\mathcal X}_k,\overset{\circ}{\mathcal X}_k,\xi_k)$,
\[
\begin{aligned}
V_{ki}: {\mathcal X}_k \otimes {\mathcal Y}_i(l,k) \otimes {\mathcal Y}(l,i) &\to {\mathcal X} \\
W_{ik}: {\mathcal Y}(r,i) \otimes {\mathcal Y}_i(r,k) \otimes {\mathcal X}_k &\to {\mathcal X}
\end{aligned}
\]
if $k \in K_i$ (note that ${\mathcal X} = {\mathcal Y}$ here). Moreover we have that
\[
\begin{aligned}
\sigma\lambda_i(\lambda_k(T)) &= \lambda_k(T)\sigma \\
\sigma\rho_i(\rho_k(T)) &= \rho_k(T)\sigma
\end{aligned}
\]
where $T \in {\mathcal L}({\mathcal X}_k)$, $k \in K_i$ and the $\lambda_k,\rho_k$ in the left-hand sides are operators on ${\mathcal Y}_i$, while the $\lambda_k,\rho_k$ in the right-hand sides are operators on ${\mathcal X}$.

\section{Bi-freeness for Pairs of Faces}
\label{sec2}

\noindent
{\bf 2.1. Definition.} A pair of faces (or face-pair) in a non-commutative probability space $({\mathcal A},\varphi)$ is an ordered pair $(({\mathcal B},\beta),({\mathcal C},\gamma))$ where ${\mathcal B},{\mathcal C}$ are unital algebras and $\beta: {\mathcal B} \to {\mathcal A}$, $\gamma: {\mathcal C} \to {\mathcal A}$ are unital homomorphisms. We shall refer to $({\mathcal B},\beta)$ and to $({\mathcal C},\gamma)$ as the left face and respectively as the right face in the pair of faces. If $1 \in {\mathcal B}$ and $1 \in {\mathcal C}$ are subalgebras in ${\mathcal A}$ and $\beta,\gamma$ are the inclusion homomorphisms, the pair $(({\mathcal B},\beta),({\mathcal C},\gamma))$ will also be denoted $({\mathcal B},{\mathcal C})$ and ${\mathcal B},{\mathcal C}$ will be called the right and left faces.

\bigskip
\noindent
{\bf 2.2. Definition.} A two-faced family of non-commutative random variables in a non-commutative probability space $({\mathcal A},\varphi)$ is an ordered pair $((b_i)_{i \in I},(c_j)_{j \in J})$ of families of non-commutative random variables in $({\mathcal A},\varphi)$ (i.e., the $b_i$ and $c_j$ are elements of ${\mathcal A}$). The associated pair of faces is $(({\mathbb C}\langle X_i \mid i \in I\rangle,\beta),({\mathbb C}\langle Y_j \mid j \in J\rangle,\gamma))$ where $\beta,\gamma$ are the unital homomorphisms so that $\beta(X_i) = b_i$, $\gamma(Y_j) = c_j$, $i \in I$, $j \in J$. We shall refer to $(b_i)_{i \in I}$ as the family of left variables and $(c_j)_{j \in J}$ as the family of right variables in the two-faced family.

\bigskip
\noindent
{\bf 2.3.} We pass to the definitions for distributions and {\em we shall use free products of unital algebras, which will always be meant to be with amalgamation over ${\mathbb C}1$}. The free product of family $({\mathcal A}_i)_{i \in I}$ of unital algebras will be denoted by $\underset{i \in I}{*} {\mathcal A}_i$.

\bigskip
\noindent
{\bf 2.4. Definition.} If $((({\mathcal B}_k,\beta_k),({\mathcal C}_k,\gamma_k)))_{k \in K} = \pi$ is a family of pairs of faces in $({\mathcal A},\varphi)$, then is joint distribution is the functional $\mu_{\pi}: \underset{k \in K}{*} ({\mathcal B}_k * {\mathcal C}_k) \to {\mathbb C}$ defined by $\mu_{\pi} = \varphi \circ \alpha$, where $\alpha: \underset{k \in K}{*} ({\mathcal B}_k * {\mathcal C}_k) \to {\mathcal A}$ is the homomorphism such that $\alpha|_{{\mathcal B}_k} = \beta_k$, $\alpha|_{{\mathcal C}_k} = \gamma_k$. The distribution of a pair of faces is defined as the distribution of the singleton family of pairs it represents.

\bigskip
\noindent
{\bf 2.5. Definition.} If $((b_i)_{i \in I},(c_j)_{j \in J}) = ({\hat b},{\hat c})$ is a two-faced family of non-commutative random variables in $({\mathcal A},\varphi)$, the its distribution $\mu_{{\hat b},{\hat c}}$ is the functional
\[
\mu_{{\hat b},{\hat c}}: {\mathbb C}\langle X_i,Y_j \mid i \in I,\ j \in J\rangle \to {\mathbb C},
\]
defined as $\mu_{{\hat b},{\hat c}} = \varphi \circ \alpha$, where $\alpha: {\mathbb C}\langle X_i,Y_j \mid i \in I,\ j \in J\rangle \to {\mathcal A}$ is the homomorphism so that $\alpha(X_i) = b_i$, $\alpha(Y_j) = c_j$.

\bigskip
\noindent
{\bf 2.6. Definition.} A family $\pi = ((({\mathcal B}_k,\beta_k),({\mathcal C}_k,\gamma_k)))_{k \in K}$ of pairs of faces in $({\mathcal A},\varphi)$ is {\em bi-freely independent} (abbreviated {\em bi-free}) if there is a family of vector spaces with specified state vector $(({\mathcal X}_k,\overset{\circ}{\mathcal X}_k,\xi_k))_{k \in K}$ and unital homomorphisms $l_k: {\mathcal B}_k \to {\mathcal L}({\mathcal X}_k)$, $r_k: {\mathcal C}_k \to {\mathcal L}({\mathcal X}_k)$, such that if ${\tilde \pi} = ((({\mathcal B}_k,\lambda_k \circ l_k),({\mathcal C}_k,\rho_k \circ r_k)))_{k \in K}$ is the family of faces of $({\mathcal L}({\mathcal X}),p)$ where $({\mathcal X},\overset{\circ}{\mathcal X},\xi) = \underset{k \in K}{*} ({\mathcal X}_k,\overset{\circ}{\mathcal X}_k,\xi_k)$, we have the equality of distributions $\mu_{\pi} = \mu_{\tilde \pi}$. The family $\pi$ is {\em strictly bi-free} if it is bi-free and $k \ne k' \Rightarrow [\beta_k({\mathcal B}_k),\gamma_{k'},({\mathcal C}_{k'})] = 0$.

\bigskip
\noindent
{\bf 2.7.} We will need to show that in a certain sense the preceding definition does not depend on a particular choice of $l_k,r_k$ and $({\mathcal X}_k,\overset{\circ}{\mathcal X}_k,\xi_k)$. We shall use $1.11$, $1.12$, $1.13$ in dealing with this matter. We begin with the following lemma.

\bigskip
\noindent
{\bf 2.8. Lemma.} {\em Let ${\mathcal B}_k,{\mathcal C}_k$ be unital algebras and let $l_k: {\mathcal B}_k \to {\mathcal L}({\mathcal X}_k)$, $r_k: {\mathcal C}_k \to {\mathcal L}({\mathcal X}_k)$, $l'_k: {\mathcal B}_k \to {\mathcal L}({\mathcal X}'_k)$, $r'_k: {\mathcal C} \to {\mathcal L}({\mathcal X}'_k)$ be unital homomorphisms, where $({\mathcal X}_k,\overset{\circ}{\mathcal X}_k,\xi_k)$, $({\mathcal X}'_k,\overset{\circ}{\mathcal X}{}'_k,\xi'_k)$ are vector spaces with specified state-vectors and assume there are morphisms $S_k: ({\mathcal X}_k,\overset{\circ}{\mathcal X}_k,\xi_k) \to ({\mathcal X}'_k,\overset{\circ}{\mathcal X}{}'_k,\xi'_k)$ so that $S_k$ intertwines the representations $l_k$ and $l'_k$ of ${\mathcal B}_k$ and also the representations $r_k,r'_k$ of ${\mathcal C}_k$. Then, if $({\mathcal X},\overset{\circ}{\mathcal X},\xi) = \underset{k \in K}{*} ({\mathcal X}_k,\overset{\circ}{\mathcal X}_k,\xi_k)$, $({\mathcal X}',\overset{\circ}{\mathcal X}{}',\xi') = \underset{k \in K}{*} ({\mathcal X}'_k,\overset{\circ}{\mathcal X}{}'_k,\xi'_k)$, $S = \underset{k \in K}{*} S_k$ and $\alpha: \underset{k \in K}{*} ({\mathcal B}_k * {\mathcal C}_k) \to {\mathcal L}({\mathcal X})$, $\alpha': \underset{k \in K}{*} ({\mathcal B}_k * {\mathcal C}_k) \to {\mathcal L}({\mathcal X}')$ are the homomorphisms so that $\alpha|_{{\mathcal B}_k} = \lambda_k \circ l_k$, $\alpha|_{{\mathcal C}_k} = \rho_k \circ r_k$, $\alpha'|_{{\mathcal B}_k} = \lambda'_k \circ l'_k$, $\alpha'|_{{\mathcal C}_k} = \rho'_k \circ r'_k$ it follows that $S$ intertwines $\alpha$ and $\alpha'$. Moreover, we also have $\varphi_p \circ \alpha = \varphi_{p'} \circ \alpha'$ where $\varphi_p,\varphi_{p'}$, are the expectation functionals on ${\mathcal L}({\mathcal X})$ and ${\mathcal L}({\mathcal X}')$.}

\bigskip
\noindent
{\bf {\em Proof.}} If $b \in {\mathcal B}_k$ then $S_kl_k(b) = l'_k(b)S_k$, by the intertwining assumption and it follows from Lemma~$1.13$ that $S\lambda_k(l_k(b)) = \lambda'_k(l'_k(b))S$. Similarly we show that $S\rho_k(r_k(c)) = \rho'_k(r'_k(c))S$ if $c \in {\mathcal C}_k$. This means that $S\alpha(x) = \alpha'(x)S$ when $x$ runs over the sets ${\mathcal B}_k$, ${\mathcal C}_k$, $k \in K$, which generate $\underset{k \in K}{*} ({\mathcal B}_k * {\mathcal C}_k)$ and hence $S$ intertwines the representations $\alpha$ and $\alpha'$ of $\underset{k \in K}{*} ({\mathcal B}_k * {\mathcal C}_k)$ on ${\mathcal X}$ and ${\mathcal X}'$. Since $S: ({\mathcal X},\overset{\circ}{\mathcal X},\xi) \to ({\mathcal X}',\overset{\circ}{\mathcal X}{}',\xi')$ is a morphism the fact that $\varphi_p \circ \alpha = \varphi_{p'} \circ \alpha'$ is then a consequence of Lemma~$1.11$.\qed

\bigskip
\noindent
{\bf 2.9. Proposition.} {\em Let $\pi = ((\pi_k))_{k \in K} = ((({\mathcal B}_k,\beta_k),({\mathcal C}_k,\gamma_k)))_{k \in K}$ be a family of pairs of faces in $({\mathcal A},\varphi)$ and for each $k \in K$ let $({\mathcal X}'_k,\overset{\circ}{\mathcal X}{}'_k,\xi'_k)$ be a vector space with specified state-vector and let $l'_k: {\mathcal B}_k \to {\mathcal L}({\mathcal X}'_k)$, $r'_k: {\mathcal C}_k \to {\mathcal L}({\mathcal X}'_k)$ be unital homomorphisms which define a pair of faces $\pi'_k: (({\mathcal B}_k,l'_k),({\mathcal C}_k,r'_k)$ in $({\mathcal L}({\mathcal X}'_k),\varphi_{p'_k})$ so that we have the equalities of distributions $\mu_{\pi_k} = \mu_{\pi'_k}$, $k \in K$. Let further $({\mathcal X}',\overset{\circ}{\mathcal X}{}',\xi') = \underset{k \in K}{*} ({\mathcal X}'_k,\overset{\circ}{\mathcal X}{}'_k,\xi'_k)$ and let ${\tilde \pi}' = ({\tilde \pi}'_k)_{k \in K} = ((({\mathcal B}_k,\lambda'_k \circ l'_k),({\mathcal C}_k,\rho'_k \circ r'_k)))_{k \in K}$ be the family of pairs of faces in $({\mathcal L}({\mathcal X}'),\varphi_{p'})$. Then $\pi$ is bi-free iff $\mu_{\pi} = \mu_{\tilde \pi}$.}

\bigskip
\noindent
{\bf {\em Proof.}} That $\pi$ is bi-free if $\mu_{\pi} = \mu_{\tilde \pi}$, is precisely the definition of bi-freeness. The converse amounts to the fact that if $({\mathcal X}''_k,\overset{\circ}{\mathcal X}{}''_k,\xi''_k)$, $l''_k: {\mathcal B}_k \to {\mathcal L}({\mathcal X}''_k)$, $r''_k: {\mathcal C}_k \to {\mathcal L}({\mathcal X}''_k)$ are so that $\pi''_k = (({\mathcal B}_k,l''_k),({\mathcal C}_k,r''_k))$ in $({\mathcal L}({\mathcal X}''_k),\varphi_{p''_k})$ satisfies $\mu_{\pi''_k} = \mu_{\pi'_k}$ for each $k \in K$, then we must have $\mu_{{\tilde \pi}'} = \mu_{{\tilde \mu}''}$ where ${\tilde \pi}'' = ({\tilde \pi}''_k)_{k \in K} = ((({\mathcal B}_k,\lambda''_k \circ l''_l),({\mathcal C}_k,\rho''_k \circ r''_k)))_{k \in K}$ in $({\mathcal L}({\mathcal X}''),\varphi_{p''})$ and $({\mathcal X}'',\overset{\circ}{\mathcal X}{}'',\xi'') = \underset{k \in K}{*} ({\mathcal X}''_k,\overset{\circ}{\mathcal X}{}''_k,\xi''_k)$. Lemma~$2.8$ provides an affirmative answer if there are morphisms
\[
S_k: ({\mathcal X}'_k,\overset{\circ}{\mathcal X}{}'_k,\xi'_k) \to ({\mathcal X}''_k,\overset{\circ}{\mathcal X}{}''_k,\xi''_k)
\]
so that $S_k$ intertwines $l'_k$ and $l''_k$ and also intertwines $r'_k$ and $r''_k$ for each $k \in K$. Of course, we cannot expect to have intertwining morphisms for every pair $(({\mathcal X}'_k,\overset{\circ}{\mathcal X}{}'_k,\xi'_k))_{k \in K}$, $(({\mathcal X}''_k,\overset{\circ}{\mathcal X}{}''_k,\xi'')_{k \in K}$ however for the desired conclusion it suffices to have some ``universal'' $(({\mathcal X}_k,\overset{\circ}{\mathcal X}_k,\xi_k))_{k \in K}$ and $l_k,r_k,k \in K$ for which there are intertwining morphisms $S_k$ to all $({\mathcal X}'_k,\overset{\circ}{\mathcal X}{}'_k,\xi'_k)$. Indeed we can take ${\mathcal X}_k = {\mathcal B}_k * {\mathcal C}_k$, $\overset{\circ}{\mathcal X}_k = \ker \mu_{\pi_k}$, $\xi_k = 1$ and define $l_k,r_k$ by the left actions of ${\mathcal B}_k$ and respectively ${\mathcal C}_k$ on ${\mathcal B}_k * {\mathcal C}_k$. We then define $S_k$ by $S_k(x) = \alpha'_k(x)\xi'_k$ where $\alpha'_k: {\mathcal B}_k * {\mathcal C}_k \to {\mathcal L}({\mathcal X}'_k)$ is the unital homomorphism so that $\alpha'_k|_{{\mathcal B}_k} = l'_k$, $\alpha'_k|_{{\mathcal C}_k} = r'_k$. If $x \in \overset{\circ}{\mathcal X}_k = \ker \mu_{\pi_k} = \ker \mu_{\pi'_k}$ we have $0 = \mu_{\pi'_k}(x)p'_k = \varphi_{p'_k}(\alpha'_k(x))p'_k = p'_k\alpha'_k(x)p'_k$ and hence $0 = p'_k\alpha'_k(x)\xi'_k$, that is $S_k\overset{\circ}{\mathcal X}_k \subset \overset{\circ}{\mathcal X}{}'_k$. It is also obvious that $S_k1 = \xi'_k$ and $S_k$ intertwines $l_k,l'_k$ and also $r_k,r'_k$.\qed

\bigskip
\noindent
{\bf 2.10. Corollary.} a) {\em If $\pi_k = (({\mathcal B}_k,\beta_k),({\mathcal C}_k,\gamma_k))$ are pairs of faces in $({\mathcal A},\varphi)$ for all $k \in K$, and the family $(\pi_k)_{k \in K} = \pi$ is bi-free, then its distribution $\mu_{\pi}$ is completely determined by the distributions $\mu_{\pi_k}$, $k \in K$.}

b) {\em Let $({\mathcal B}_k,{\mathcal C}_k)$ be pairs of unital algebras, $k \in K$ and $\mu_k: {\mathcal B}_k * {\mathcal C}_k \to {\mathbb C}$ linear functionals so that $\mu_k(1) = 1$. Then there exists a unique linear map $\mu: \underset{k \in K}{*} ({\mathcal B}_k * {\mathcal C}_k) \to {\mathbb C}$, $\mu(1) = 1$ so that the family of pairs of faces $(({\mathcal B}_k,{\mathcal C}_k))_{k \in K}$ is bi-free in $(\underset{k \in K}{*} ({\mathcal B}_k * {\mathcal C}_k),\mu)$ and $\mu|_{{\mathcal B}_k * {\mathcal C}_k} = \mu_k$ for all $k \in K$.}

\bigskip
\noindent
{\bf 2.11.} If ${\mathcal A}$ is unital algebra its algebraic state space $\Sigma({\mathcal A})$ is the convex subset of the dual of ${\mathcal A}$, consisting of linear maps $\varphi: {\mathcal A} \to {\mathbb C}$ so that $\varphi(1) = 1$. The space of distributions for families of pairs of faces arising from a family of pairs of unital algebras $(({\mathcal B}_k,{\mathcal C}_k))_{k \in K}$ is the set $\Sigma(({\mathcal B}_k,{\mathcal C}_k) \mid k \in K) = \Sigma(\underset{k \in K}{*} ({\mathcal B}_k * {\mathcal C}_k))$. For every family $\mu_k \in \Sigma({\mathcal B}_k,{\mathcal C}_k)$, $k \in K$ there is a unique {\em bi-free product distribution} $\mu = \underset{k \in K}{**} \mu_k \in \Sigma(({\mathcal B}_k,{\mathcal C}_k) \mid k \in K)$ so that $\mu|_{{\mathcal B}_k * {\mathcal C}_k} = \mu_k$, $k \in K$ and the family $(({\mathcal B}_k,{\mathcal C}_k))_{k \in K}$ is bi-free in $(\underset{k \in K}{*} ({\mathcal B}_k * {\mathcal C}_k),\mu)$ as a consequence of Corollary~$2.10$b).

\bigskip
\noindent
{\bf 2.12.} The family $((({\mathcal B}_k,\beta_k),({\mathcal C}_k,\gamma_k)))_{k \in K}$ is bi-free in \linebreak $({\mathcal A},\varphi)$ iff $((\beta_k({\mathcal B}_k),\gamma_k({\mathcal C}_k)))_{k \in K}$ is bi-free in $({\mathcal A},\varphi)$. The ``if'' is obvious, while the ``only if'' uses Proposition~$2.9$. Indeed, by Proposition~$2.9$ we can choose the vector spaces with specified state-vector for the bi-freeness of $((({\mathcal B}_k,\beta_k),({\mathcal C}_k,\gamma_k)))_{k \in K}$ to be ${\mathcal X}_k = \beta_k({\mathcal B}_k) * \gamma_k({\mathcal C}_k)$, $\xi_k = 1$ and $\overset{\circ}{\mathcal X}_k = \ker \mu_k$, where $\mu_k$ is the distribution of $(\beta_k({\mathcal B}_k),\gamma_k({\mathcal C}_k))$ in $({\mathcal A},\varphi)$ The actions of ${\mathcal B}_k$ and ${\mathcal C}_k$ on ${\mathcal X}_k$ factor through the left actions of $\beta_k({\mathcal B}_k)$ and $\gamma_k({\mathcal C}_k)$. Further, the homomorphisms ${\mathcal B}_k * {\mathcal C}_k \to \beta_k({\mathcal B}_k) * \gamma_k({\mathcal C}_k)$ then provide suitable morphisms with the intertwining properties of the actions. We leave the details as an exercise for the reader.

\bigskip
\noindent
{\bf 2.13. Proposition.} {\em Let $\pi = (({\mathcal B}_k,\beta_k),({\mathcal C}_k,\gamma_k))_{k \in K}$ be a family of pairs of faces in $({\mathcal A},\varphi)$, let $K = \underset{i \in I}{\coprod} K_i$, ${\mathcal D}_i = \underset{k \in K_i}{*} {\mathcal B}_k$, ${\mathcal E}_i = \underset{k \in K_i}{*} {\mathcal C}_k$ and let $\delta_i: {\mathcal D}_i \to {\mathcal A}$, $\epsilon_i : {\mathcal E}_i \to {\mathcal A}$ be the unital homomorphisms such that $\delta_i|_{{\mathcal B}_k} = \beta_k$, $\epsilon_i|_{{\mathcal C}_k} = \gamma_k$ if $k \in K_i$. If $\pi$ is bi-free, then $(({\mathcal D}_i,\delta_i),({\mathcal E}_i,\epsilon_i))_{i \in I}$ is bi-free.}

\bigskip
The proof is derived from Remark~$1.13$ and will be omitted.

\bigskip
\noindent
{\bf 2.14. Proposition.} {\em Let $K = \underset{i \in I}{\coprod} K_i$ and let $1 \in {\mathcal B}_k$, $1 \in {\mathcal C}_k$, $1 \in {\mathcal D}_i$, $1 \in {\mathcal E}_i$ be subalgebras in $({\mathcal A},\varphi)$ so that ${\mathcal B}_k \in {\mathcal D}_i$, ${\mathcal C}_k \in {\mathcal E}_i$ if $k \in K_i$. Assume for each $i \in I$ that $(({\mathcal B}_k,{\mathcal C}_k))_{k \in K_i}$ is bi-free in $({\mathcal D}_i \vee {\mathcal E}_i,\varphi \mid {\mathcal D}_i \vee {\mathcal E}_i)$ and that $(({\mathcal D}_i,{\mathcal E}_i))_{i \in I}$ is bi-free in $({\mathcal A},\varphi)$. Then $(({\mathcal B}_k,{\mathcal C}_k))_{k \in K}$ is bi-free in $({\mathcal A},\varphi)$.}

\bigskip
The proof of this proposition is also derived from Remark~$1.13$ and will be omitted.

\bigskip
\noindent
{\bf 2.15. Proposition.} {\em Let $1 \in {\mathcal B}_k$, $1 \in {\mathcal C}_k$, $k \in K$ be subalgebras in $({\mathcal A},\varphi)$.}

a) {\em Assume $(({\mathcal B}_k,{\mathcal C}_k))_{k \in K}$ is bi-free in $({\mathcal A},\varphi)$. Then $({\mathcal B}_k)_{k \in K}$ is free in $({\mathcal A},\varphi)$. Similarly also $({\mathcal C}_k)_{k \in K}$ is free in $({\mathcal A},\varphi)$.}

b) {\em Assume $({\mathcal B}_k)_{k \in K}$ is free in $({\mathcal A},\varphi)$ and assume also $({\mathcal C}_k)_{k \in K}$ is free in $({\mathcal A},\varphi)$. Then the families $(({\mathcal B}_k,{\mathbb C}1))_{k \in K}$ and $(({\mathbb C}1,{\mathcal C}_k))_{k \in K}$ are both bi-free in $({\mathcal A},\varphi)$.}

\bigskip
The proof is an easy free probability exercise and will be omitted.

\bigskip
\noindent
{\bf 2.16. Proposition.} a) {\em If $(({\mathcal B}_k,{\mathcal C}_k))_{k \in K}$ is bi-free in $({\mathcal A},\varphi)$ and if $L \subset K$ then the subalgebras $\underset{k \in L}{\bigvee} {\mathcal B}_k$ and $\underset{k \in K\backslash L}{\bigvee} {\mathcal C}_k$ are classically independent in $({\mathcal A},\varphi)$.}

b) {\em If $1 \in {\mathcal B}$, $1 \in {\mathcal C}$ are classically independent subalgebras in $({\mathcal A},\varphi)$ then $(({\mathcal B},{\mathbb C}1),({\mathbb C}1,{\mathcal C}))$ is bi-free in $({\mathcal A},\varphi)$.}

\bigskip
\noindent
{\bf {\em Proof.}} a) In view of $2.12$ and of Proposition~$2.13$ the family consisting of the two pairs of faces $\left( \underset{k \in K}{\bigvee} {\mathcal B}_k, \underset{k \in K}{\bigvee} {\mathcal C}_k\right)$, $\left( \underset{k \in K\backslash L}{\bigvee} {\mathcal B}_k, \underset{k \in K\backslash L}{\bigvee} {\mathcal C}_k\right)$ is bi-free in $({\mathcal A},\varphi)$ and the proof reduces to the case when $K = \{1,2\}$ and $L = \{1\}$. The assertion follows then from the following remark. Let $({\mathcal X},\overset{\circ}{\mathcal X},\xi) = \underset{k \in \{1,2\}}{*} ({\mathcal X}_k,\overset{\circ}{\mathcal X}_k,\xi_k)$ and $T_k \in {\mathcal L}({\mathcal X}_k)$, $\overset{\circ}{T}_k = {T}_k - \varphi_{p_k}({T}_k)I$, $t_k = \varphi_{p_k}({T}_k)$, $\overset{\circ}{x}_k = \overset{\circ}{T}_k\xi_k \in \overset{\circ}{\mathcal X}_k$. We have
\[
\begin{aligned}
\lambda_1({T}_1)\rho_2({T}_2)\xi &= \lambda_1({T}_1)(\overset{\circ}{x}_2 + t_2\xi) \\
&= \overset{\circ}{x}_1 \otimes \overset{\circ}{x}_2 + t_2\overset{\circ}{x}_1 + t_1\overset{\circ}{x}_2 + t_1t_2\xi
\end{aligned}
\]
so that $\varphi_p(\lambda_1({T}_1)\rho_2({T}_2)) = t_1t_2 = \varphi_{p_1}({T}_1)\varphi_{p_2}({T}_2) = \varphi_p(\lambda_1({T}_1))\varphi_p(\rho_2({T}_2))$.

b) This also essentially follows from the computation in the proof of a) showing that $\lambda_1({T}_1),\rho_2({T}_2)$ are classically independent in $({\mathcal L}({\mathcal X}),\varphi_p)$.\qed

\bigskip
\noindent
{\bf 2.17. Remark.} The notion of bi-freeness for families of two-faced families of non-commutative random variables amounts to the bi-freeness of the family of associated pairs of faces (Definition~$2.2$). By passing to the associated pairs of faces the results of this section translate into results for two-faced families of non-commutative random variables. In this setting it is also useful to describe distributions in terms of {\em two-faced moments}. For a distribution $\mu_{{\hat b},{\hat c}}$ (see Definition~$2.5$) these are the numbers $\mu_{{\hat b},{\hat c}}(M)$, whee $M$ runs through the non-commutative monomials in ${\mathbb C}\langle X_i,Y_j \mid i \in I,\ j \in J\rangle$.

\bigskip
We shall conclude the section, pointing out that like freeness, also bi-freeness amounts to algebraic relations among expectation values.

\bigskip
\noindent
{\bf 2.18. Proposition.} {\em Given sets $I,J,{T}$ for each $n \in {\mathbb N}$ and pair of $(\alpha,t): \{1,\dots,n\} \to K \times {T}$, where $K = I \coprod J$ there is $m \in {\mathbb N}$ maps $\alpha_r: \{1,\dots,n(r)\} \to K$, $t(r) \in {T}$, $1 \le r \le m$ and a universal polynomial ${\mathcal P}_{\alpha,t}(X_1,\dots,X_m)$ ($m,n(r),\alpha_r,t(r)$ depend on $(\alpha,t)$) so that:

a family of two-faced families $z_t = ((z_{t,i})_{i \in I},(z_{t,j})_{j \in J})$, $t \in {T}$ in $({\mathcal A},\varphi)$ is bi-free

iff

\noindent
for each $(\alpha,t)$ as above we have
\[
\varphi(z_{t(1),\alpha(1)}\dots z_{t(n),\alpha(n)}) = {\mathcal P}_{\alpha,t}(\varphi(z_{t(r),\alpha_r(1)} \dots z_{t(r),\alpha_r(n_r)})(1 \le r \le m).
\]
}

\bigskip
\noindent
{\bf {\em Proof.}} If $(z_t)_{t \in {T}}$ is bi-free, there is a realization as left and right variables on a free product of vector spaces. Applying recurrently $1.9$ to the computation of $z_{t(1),\alpha(1)}\dots z_{t(n),\alpha(n)}\xi$ in this realization one finds that $\varphi(z_{t(1),\alpha(1)}\dots z_{t(n),\alpha(n)})$ which is the coefficient of the specified vector $\xi$ in $z_{t(1),\alpha(1)}\dots z_{t(n),\alpha(n)}\xi$ is indeed given by a universal polynomial in moments of the $z_t$'s. We leave the details of the proof to the reader, especially since in the case $|{T}| = 2$, which does not essentially differ from the general case, a more precise result is established in detail in Lemma~$5.2$. This shows the ``only if'' part.

On the other hand that the condition is sufficient for bi-freeness, follows immediately since bi-freeness by definition means the joint distribution of $(z_t)_{t \in {T}}$ is equal to a certain distribution completely determined by the distributions of the individual $z_t$'s, so it is the one computed using the polynomials ${\mathcal P}_{\alpha,t}$.\qed

\section{Bi-freeness in $C^*$-probability Spaces}
\label{sec3}

\noindent
{\bf 3.1.} After introducing bi-freeness in the purely algebraic framework of the preceding section, we pass now to $C^*$-algebras, where there is positivity and where bi-freeness is also one of the natural probabilistic independence notions. If $({\mathcal A},\varphi)$ is a $C^*$-probability space, that is, ${\mathcal A}$ is a unital $C^*$-algebra and $\varphi$ is a state, a family of pairs of faces $(({\mathcal B}_k,\beta_k)({\mathcal C}_k,\gamma_k))_{k \in K}$ is a {\em family of pairs of $C^*$-faces} if ${\mathcal B}_k,{\mathcal C}_k$ are $C^*$-algebras and $\beta_k,\gamma_k$ are unital $*$-homomorphism. Similarly we have corresponding $W^*$-notions, where $({\mathcal A},\varphi)$ is a $W^*$-probability space, the ${\mathcal B}_k,{\mathcal C}_k$ are $W^*$-algebras and the unital $*$-homomorphisms $\beta_k,\gamma_k$ are ultraweakly continuous. If $\pi$ is a family of pairs of $C^*$-faces in a $C^*$-probability space then {\em the distribution $\mu_{\pi}$ extends by continuity to the full free-$C^*$-algebra product and is a state on this $C^*$-algebra}. Also in the $C^*$-algebra setting the two-faced families of non-commutative random variables often have specifications that some of variables be hermitian, unitary, etc.

\bigskip
\noindent
{\bf 3.2.} If $\pi = (({\mathcal B}_k,\beta_k),({\mathcal C}_k,\gamma_i))_{k \in K}$ is a family of pairs of $C^*$-faces in a $C^*$-probability space $({\mathcal A},\varphi)$, then the free product of vector spaces with specified state-vector which are used in the bi-freeness condition, can be chosen to be a free product of Hilbert spaces with specified unit vector. More precisely, $\varphi \circ \beta_k$ and $\varphi \circ \gamma_k$ are states on ${\mathcal B}_k$ and ${\mathcal C}_k$ and therefore $(\varphi \circ \beta_k)*(\varphi \circ \gamma_k) = \psi_k$ is a state o the full free product $C^*$-algebra ${\mathcal B}_k * {\mathcal C}_k$ (with amalgamation over ${\mathbb C} 1$). The GNS-construction applied to ${\mathcal B}_k * {\mathcal C}_k$ and $\psi_k$ gives rise to a $*$-representation of ${\mathcal B}_k * {\mathcal C}_k$ on a Hilbert space ${\mathcal H}_k$ with specified unit vector $\xi_k$, the restrictions of which to ${\mathcal B}_k$ and ${\mathcal C}_k$ give rise to $*$-representations $l_k,r_k$ of these algebras on ${\mathcal H}_k$. Since $\lambda_k$ and $\rho_k$ in the context of $({\mathcal H},\xi) = {\underset{k \in K}{*}} ({\mathcal H}_k,\xi_k)$ are $*$-representations, we have $*$-representations $\lambda_k \circ l_k$, $\rho_k \circ r_k$ of ${\mathcal B}_k$ and ${\mathcal C}_k$ on $({\mathcal H},\xi)$. It follows that the joint distribution $\mu_{\pi}$ is the pull-back of the vector state $\langle \cdot\xi,\xi\rangle$ via a $*$-representation of ${\underset{k \in K}{*}} ({\mathcal B}_k * {\mathcal C}_k)$ on ${\mathcal H}$. This shows that {\em the distribution $\mu_{\pi}$ is a state on the $C^*$-algebra ${\underset{k \in K}{*}} ({\mathcal B}_k * {\mathcal C}_k)$} (here the free products are full free products of unital $C^*$-algebras). Thus {\em in the $C^*$-context the bi-free product of distributions ${\underset{k \in K}{**}} \mu_k$ extends by continuity to a state on the $C^*$-algebra free product}. This is the same as saying that {\em if $\mu_k$ are states of ${\mathcal B}_k * {\mathcal C}_k$ there is an unique bi-free product state $\mu$ on ${\underset{k \in K}{*}} ({\mathcal B}_k * {\mathcal C}_k)$ so that the family $(({\mathcal B}_k,{\mathcal C}_k))_{k \in K}$ is bi-free. The bi-free product state will also be denoted by ${\underset{k \in K}{**}} \mu_k$, like its algebraic relative}.

\bigskip
\noindent
{\bf 3.3. Proposition.} a) {\em If the family of two-faced pairs of algebras $(({\mathcal B}_k,{\mathcal C}_k))_{k \in K}$ in the Banach probability space $({\mathcal A},\varphi)$ is bi-free then the family of norm-closures $((\overline{{\mathcal B}}_k,\overline{{\mathcal C}}_k))_{k \in k}$ is bi-free in $({\mathcal A},\varphi)$.}

b) {\em If the family of two-faced pairs of $*$-subalgebras $(({\mathcal B}_k,{\mathcal C}_k))_{k \in K}$ in the $C^*$-probability space $({\mathcal A},\varphi)$ is bi-free, then the family of pairs of $C^*$-algebras $((\overline{{\mathcal B}}_k,\overline{{\mathcal C}}_k))_{k \in K}$ is bi-free.}

c) {\em If the family of two-faced pairs of $*$-subalgebras $(({\mathcal B}_k,{\mathcal C}_k))_{k \in K}$ is bi-free in the $W^*$-probability space $({\mathcal A},\varphi)$ is bi-free, the family of pairs $((\overline{{\mathcal B}}_k^w,\overline{{\mathcal C}}_k^w))_{k \in K}$ is also bi-free in $({\mathcal A},\varphi)$.}

\bigskip
Like for freeness, bi-freeness reduces to the bi-freeness of families of left and right variables in the algebras and this amounts to polynomial relations among moments (Proposition~$2.18$). The polynomial relations among moments are conserved when taking norm-limits of the variables and also when taking $*$-ultrastrong limits of bounded nets after using the Kaplansky density theorem for c).

\newpage
\section{Bi-free Convolutions}
\label{sec4}

\noindent
{\bf 4.1. Definition.} If the pair of two-faced families of non-commutative random variables $({\hat b},{\hat c}) = ((b_i)_{i \in I},(c_j)_{j \in J}$, $({\hat d},{\hat e}) = ((d_i)_{i \in I},(e_j)_{j \in J})$ is bi-free in $({\mathcal A},\varphi)$ then their joint distribution being completely determined by the distributions $\mu_{{\hat b},{\hat c}}$ and $\mu_{{\hat d},{\hat e}}$, in particular the distributions $\mu_{{\hat b}+{\hat d},{\hat c}+{\hat e}},\mu_{{\hat b}{\hat d},{\hat c}{\hat e}}$ of $((b_i+d_i)_{i \in I},(c_j+e_j)_{j \in J})$ and $((b_id_i)_{i \in I},(c_je_j)_{j \in J})$ are completely determined by $\mu_{{\hat b},{\hat c}},\mu_{{\hat d},{\hat e}}$. This defines additive and multiplicative bi-free convolution operations $\boxplus\boxplus$ and $\boxtimes\boxtimes$ on distributions of two-faced families of non-commutative random variables with pair of index sets $(I,J)$ so that
\[
\begin{aligned}
\mu_{{\hat b}+{\hat d},{\hat c}+{\hat e}} &= \mu_{{\hat b},{\hat c}} \boxplus\boxplus \mu_{{\hat d},{\hat e}} \\
\mu_{{\hat b}{\hat d}+{\hat c}{\hat e}} U&= \mu_{{\hat b},{\hat c}} \boxtimes\boxtimes \mu_{{\hat d},{\hat e}}.
\end{aligned}
\]

\bigskip
\noindent
{\bf 4.2. Remark.} If $({\mathcal A},\varphi)$ is a $C^*$-probability space and if the $b_i,c_j,d_i,e_j$, $i \in I$, $j \in J$ are self-adjoint then also $b_i + d_i, c_j + e_j$ are self-adjoint and the joint distribution $\mu_{{\hat b}+{\hat d},{\hat c}+{\hat e}}$ is a state of the unital $*$-algebra ${\mathbb C}\langle X_i,Y_j \mid i \in I,\ j \in J\rangle$ endowed with the involution given by $X_i = X_i^*$, $Y_j = Y_j^*$, for which boundedness conditions for the variables $\mu_{{\hat b}+{\hat d},{\hat c}+{\hat e}}({\mathcal P}^*X_i^2{\mathcal P}) \le C_i\mu_{{\hat b}+{\hat d},{\hat c}+{\hat e}}({\mathcal P}^*{\mathcal P})$, $\mu_{{\hat b}+{\hat d},{\hat c}+{\hat e}}({\mathcal P}^*Y_j^2{\mathcal P}) \le C_j\mu_{{\hat b}+{\hat d},{\hat c}+{\hat e}}({\mathcal P}^*{\mathcal P})$ hold for all ${\mathcal P} \in {\mathbb C}\langle X_i,Y_j \mid i \in I,\ j \in J\rangle$. Thus if $\mu$ and $\nu$ are states of ${\mathbb C}\langle X_i,Y_j \mid i \in I,\ j \in J\rangle$ endowed with this involution and for which boundedness conditions for the variables are satisfied, then $\mu \boxplus\boxplus \nu$ is such a state, that is, the distribution of some two-faced family $((f_i))_{i \in I},(g_j)_{j \in J})$ of self-adjoint elements in a $C^*$-probability space. A convenient way to handle such distributions, in view of the varying bounds, i.e., norms for the elements, is to replace ${\mathbb C}\langle X_i,Y_j \mid i \in I,\ j \in J\rangle$ by some pro-$C^*$-algebra, like we used for free convolution in \cite{15}.

\bigskip
\noindent
{\bf 4.3. Remark.} To handle operations on $*$-distributions for two-faced families of non-self-adjoint non-commutative random variables in a $C^*$-probability space, one can pass to hermitian and antihermitian parts or proceed as follows. Instead of the two-faced family $((b_i)_{i \in I},(c_j)_{j \in J})$ we must consider the family $((b_i)_{i \in I} \coprod (b_i^*)_{i \in I},(c_j)_{j \in J} \coprod (c_j^*)_{j \in J})$ with index sets $I \coprod I,J \coprod J$. The distribution becomes a state of ${\mathbb C}\langle X_i,X_i^*,Y_j,Y_j^* \mid i \in I,\ j \in J\rangle$ with involution defined by $(X_i)^* = X_i^*$, $(X_i^*)^* = X_i$, $(Y_j)^* = Y_j^*$, $(Y_j^*)^* = Y_j$ and with appropriate boundedness conditions for the variables. Again, a pro-$C^*$-algebra framework can be adapted for this along the lines of \cite{15}.

\bigskip
\noindent
{\bf 4.4. Remark.} When dealing with multiplicative bi-free convolution of two-faced families of unitaries the result will preserve the unitary property. If we want to consider the $*$-distributions, then a convenient way for such a two-faced unitary family $((u_i)_{i \in I},(v_j)_{j \in J})$ in a $C^*$-probability space is to look at the homomorphism of the full $C^*$-algebra of the free group on generators indexed by $I \coprod J$, $C^*({\mathbb Z}^*{}^{I \coprod J})$ two which it gives rise and the $*$-distribution is a state of this $C^*$-algebra which identifies also with a certain positive-definite function on the group. Thus multiplicative bi-free convolution becomes an operation on states of $C*({\mathbb Z}^*{}^{I \coprod J})$.

\section{Bi-free Cumulants}
\label{sec5}

\noindent
{\bf 5.1.} The purpose of this section is to define and prove the existence and uniqueness of bi-free cumulants, which linearize additive bi-free convolution. It will be convenient in what follows to distinguish left and right variables only the the index sets, not by the letters we see. Thus we shall consider two-faced families of non-commutative random variables $((z_i)_{i \in I},(z_j)_{j \in J})$, where of course the sets $I$ and $J$ are presumed to be disjoint.

\bigskip
\noindent
{\bf 5.2. Lemma.} {\em Let $z' = ((z'_i)_{i \in I},(z'_j)_{j \in J})$, $z'' = ((z''_i)_{i \in I},(z''_j)_{j \in J})$ denote a bi-free pair of two-faced families of non-commutative random variables in some non-commutative probability space $({\mathcal A},\varphi)$ and let $\alpha: \{1,\dots,n\} \to I \coprod J$, $\epsilon: \{1,\dots,n\} \to \{',''\}$ be given. Then there is a universal polynomial ${\mathcal P}$ so that
\[
\begin{aligned}
&\varphi(z_{\alpha(1)}^{\epsilon(1)} \dots z_{\alpha(n)}^{\epsilon(n)}) \\
&\quad= {\mathcal P}(\varphi(z'_{\alpha(k_1)} \dots a'_{\alpha(k_r)}),\varphi(z''_{\alpha(l_1)} \dots z''_{\alpha(l_s)}) \mid 1 \le k_1 < \dots < k_r \le n, \\
&\qquad\, 1 \le l_1 < \dots < l_s \le n, \epsilon(k_a) = ',\epsilon(l_b) = '', 1 \le a \le r,1 \le b \le s).
\end{aligned}
\]
Moreover ${\mathcal P}$ has integer coefficients and if $\varphi(z'_{\alpha(k_1)} \dots z'_{\alpha(k_r)})$, $\varphi(z''_{\alpha(l_1)} \dots z''_{\alpha(l_s)})$ are associated bidegree $(r,0)$ and $(0,s)$ respectively, then ${\mathcal P}$ is homogeneous of bidegree $(p,q)$, where $p = \#\{1 \le j \le n \mid \epsilon(a) = '\}$, $q = \#\{1 \le b \le n \mid \epsilon(b) = ''\}$.}

\bigskip
The lemma will be a consequence of a more detailed lemma.

\bigskip
\noindent
{\bf 5.3. Lemma.} {\em Let $({\mathcal X},{\mathcal X}_0,\xi) = ({\mathcal X}',{\mathcal X}'_0,\xi') * ({\mathcal X}'',{\mathcal X}''_0,\xi'')$ and let $\epsilon: \{1,\dots,n\} \to \{',''\}$, $\mu: \{1,\dots,n\} \to \{\lambda,\rho\}$, $T_k \in {\mathcal L}({\mathcal X}^{\epsilon(k)})$, $1 \le k \le n$. Further, if $A \subset \{1,\dots,n\}$, $A \ne \emptyset$ and $\delta \in \{{}^{\circ},\varphi\}$, let $M(A,\delta)$ be defined as follows: if 
\[
A = \{\alpha_1 < \dots < \alpha_k,\epsilon(\alpha_1) = \dots = \epsilon(\alpha_k) = \epsilon\}
\]
then
\[
\begin{aligned}
M(A,\delta) &= T_{\alpha_1} \dots T_{\alpha_k}\xi^{\epsilon} - \varphi(T_{\alpha_1} \dots T_{\alpha_k})\xi^{\epsilon} \in {\mathcal X}_0^{\epsilon} \text{ when $\delta = 0$} \\
M(A,\delta) &= \varphi(T_{\alpha_1} \dots T_{\alpha_k}) \in {\mathbb C} \text{ when $\delta = \varphi$}
\end{aligned}
\]
and otherwise $M(A,\delta) = 0$. Then there are universal integer coefficients $c_{A_1,\dots,A_k,\delta_1,\dots,\delta_k}$, where $(A_1,\dots,A_k)$ is an ordered partition of $\{1,\dots,n\}$ into non-empty subsets, $\{\delta_1,\dots,\delta_k\} \in \{{}^{\circ},\varphi\}^k$, $1 \le k \le n$, so that
\[
\begin{aligned}
&\mu(1)_{\epsilon(1)}(T_1) \dots \mu(n)_{\epsilon(n)}(T_n)\xi \\
&\quad = \sum c_{A_1,\dots,A_k,\delta_1,\dots,\delta_k}M(A_1,\delta_1) \otimes \dots \otimes M(A_k,\delta_k)
\end{aligned}
\]
with the convention that the $\otimes$ sign is also used for multiplication with scalars or of scalars.}

\bigskip
\noindent
{\bf {\em Proof of Lemma 5.3.}} The lemma is proved by induction over $n$. The induction step from $n$ to $n+1$ is then a straightforward case-by-case application of the formulae in $1.9$ depending on whether in $\mu(1)_{\epsilon(1)}(T_1)\dots \mu(n+1)_{\epsilon(n+1)}(T_{n+1})\xi$ we have $\mu(1) = \lambda$ or $\mu(1) = \rho$. Next is in more detail.

The assumption that the Lemma holds for $n$ applied to $\{2,\dots,n+1\}$ instead of $\{1,\dots,n\}$ gives
\[
\begin{aligned}
&\mu(2)_{\epsilon(2)}(T_2)\dots \mu(n+1)_{\epsilon(n+1)}(T_{n+1})\xi \\
&\quad = \sum c_{A_1,\dots,A_k,\delta_1,\dots,\delta_k}M(A_1,\delta_1) \otimes \dots \otimes M(A_k,\delta_k)
\end{aligned}
\]
where $(A_1,\dots,A_k)$ is now a partition of $\{2,\dots,n+1\}$. We will show that $\mu(1)_{\epsilon(1)}(T_1)M(A_1,\delta_1) \otimes \dots \otimes M(A_k,\delta_k)$ is of the desired form, which will prove the induction step.

\bigskip
{\bf Case 1. $\mu(1) = \lambda$.} If $\delta_1 = \dots = \delta_k = \varphi$ all the $M$'s are scalars and this reduces to
\[
\lambda_{\epsilon(1)}(T_1)\xi = (T_1-\varphi_{\epsilon(1)}(T_1)I)\xi_1 + \varphi_{\epsilon(1)}(T_1)\xi = M(\{1\},{}^{\circ}) + M(\{1\},\varphi).
\]
In this case the partition of $\{1,\dots,n+1\}$ is $\{1\} \coprod A_1 \coprod \dots \coprod A_k$.

If not all $\delta_j = \varphi$, assume $p$ is the smallest index so that $\delta_p = {}^{\circ}$. We have
\[
\begin{aligned}
&\lambda_{\epsilon(1)}(T_1)M(A_1,\delta_1) \otimes \dots \otimes M(A_k,\varphi) \\
&\quad = M(A_1,\varphi) \otimes \dots \otimes M(A_{p-1},\varphi) \otimes (\lambda_{\epsilon(1)}(T_1)M(A_p,{}^{\circ})) \\
&\quad \otimes M(A_{p+1},\delta_{p+1}) \otimes \dots \otimes M(A_{n+1},\delta_{n+1}).
\end{aligned}
\]
If $\epsilon(1) \ne \epsilon(p)$ we have
\[
\lambda_{\epsilon(1)}(T_1)M(A_p,{}^{\circ}) = M(\{1\},{}^{\circ}) \otimes M(A_p,{}^{\circ}) + M(\{1\},\varphi) \otimes M(A_p,{}^{\circ}).
\]
If $\epsilon(1) = \epsilon(p)$ on the other hand
\[
\begin{aligned}
\lambda_{\epsilon(1)}(T_1)M(A_p,{}^{\circ}) &= T_1(T_{\alpha_1} \dots T_{\alpha_l} - \varphi(T_{\alpha_1} \dots T_{\alpha_l})I)\xi \\
&= (T_1T_{\alpha_1} \dots T_{\alpha_l} - \varphi(T_1T_{\alpha_1} \dots T_{\alpha_l})I)\xi \\
&+ \varphi(T_1T_{\alpha_1} \dots T_{\alpha_1})\xi + \varphi(T_{\alpha_1} \dots T_{\alpha_l})(T_1-\varphi(T_1)I)\xi \\
&+ \varphi(T_1)\varphi(T_{\alpha_1} \dots T_{\alpha_l}) \\
&= M(\{1\} \cup A_p,{}^{\circ}) + M(\{1\} \cup A_p,\varphi) \\
&+ M(\{1\},{}^{\circ}) \otimes M(A_p,\varphi) + M(\{1\},\varphi) \otimes M(A_p,\varphi),
\end{aligned}
\]
which concludes the proof of this case, remarking that the partitions which appear are $\{1\} \coprod A_1 \coprod \dots \coprod A_k$ and $A_1 \coprod \dots \coprod (\{1\} \cup A_p) \coprod \dots \coprod A_k$.

\bigskip
{\bf Case 2. $\mu(1) = \rho$.} Again, if $\delta_1 = \dots = \delta_k = \varphi$ this reduces to
\[
\rho_{\epsilon(1)}(T_1)\xi = M(\{1\},{}^{\circ}) + M(\{1\},\varphi).
\]

If not all $\delta_j = \varphi$, we must now look at the largest index $p$ so that $\delta_p = {}^{\circ}$. We have:
\[
\begin{aligned}
&\rho_{\epsilon(1)}(T_1)M(A_1,\delta_1) \otimes \dots \otimes M(A_k,\varphi) \\
&\quad = M(A_1,\delta_1) \otimes \dots \otimes (\rho_{\epsilon(1)}(T_1)M(A_p,{}^{\circ}) \otimes \dots \otimes M(A_{n+1},\varphi).
\end{aligned}
\]
If $\epsilon(1) \ne \epsilon(p)$ we have
\[
\rho_{\epsilon(1)}(T_1)M(A_p,{}^{\circ}) = M(A_p,{}^{\circ}) \otimes M(\{1\},{}^{\circ}) + M(A_p,{}^{\circ}) \otimes M(\{1\},\varphi).
\]
If $\epsilon(1) = \epsilon(p)$ we get
\[
\begin{aligned}
\rho_{\epsilon(1)}(T_1)M(A_p,{}^{\circ}) &= M(\{1\} \cup A_p,{}^{\circ}) + M(\{1\} \cup A_p,\varphi) \\
&+ M(\{1\},{}^{\circ}) \otimes M(A_p,\varphi) + M(\{1\},\varphi) \otimes M(A_p,\varphi)
\end{aligned}
\]
which concludes the proof. (Remark that also here we get $T_1T_{\alpha_1} \dots T_{\alpha_l}$, only the index $p$ is the right-most.)\qed

\bigskip
\noindent
{\bf 5.4. {\em Proof of Lemma 5.2.}} In view of Proposition~$2.9$ it suffices to prove the lemma when $z'_i = \lambda,(T_i)$, $z'_j = \rho,(T_j)$, $z''_i = \lambda,,(T_i)$, $z''_j = \lambda,,(T_j)$ for some $T_i,T_j$. Then $\varphi(z_{\alpha(1)}^{\epsilon(1)} \dots z_{\alpha(n)}^{\epsilon(n)})$ amounts to taking in the sum appearing in the formula for $\mu(1)_{\delta(1)}(T_1) \dots \mu(n)_{\epsilon(n)}(T_n)$ in Lemma~$5.3$ only the terms with $\delta_1 = \dots = \delta_k = \varphi$.\qed

\bigskip
\noindent
{\bf 5.5. Remark.} In view of Lemma~$5.3$ the statement of Lemma~$5.2$ can be made more precise than just saying that ${\mathcal P}$ is homogeneous of bi-degree $(p,q)$. Indeed, the products of moments $\varphi(z'_{\alpha(k_1)} \dots z'_{\alpha(k_r)})$ and $\varphi(z''_{\alpha(l_1)} \dots z''_{\alpha(l_s)})$, in view of Lemma~$5.3$, are restricted to those which correspond to an ordered partition $(A_1,\dots,A_k)$ of $\{1,\dots,n\}$.

\bigskip
\noindent
{\bf 5.6. Proposition.} {\em Let $z' = ((z'_i)_{i \in I},(z'_j)_{j \in J})$, $z'' = ((z''_i)_{i \in I},(z''_j)_{j \in J})$ denote a bi-free pair of two-faced families of non-commutative random variables in a non-commutative probability space $({\mathcal A},\varphi)$ and let $\alpha: \{1,\dots,n\} \to I \coprod J$ be given. Then there is a universal polynomial ${\tilde Q}$ so that
\[
\begin{aligned}
&\varphi((z'_{\alpha(1)} + z''_{\alpha(1)}) \dots (z'_{\alpha(n)} + z''_{\alpha(n)})) = \varphi(z'_{\alpha(1)} \dots z'_{\alpha(n)}) + \varphi(z''_{\alpha(1)} \dots z''_{\alpha(n)}) \\
&\quad + {\tilde Q}(\varphi(z'_{\alpha(k_1)} \dots z'_{\alpha(k_n)}), \\
&\quad\quad\, \varphi(z''_{\alpha(k_1)} \dots z''_{\alpha(k_r)}) \mid 1 \le k_1 < \dots < k_r \le n, 1 \le r \le n-1).
\end{aligned}
\]
More precisely ${\tilde Q}$ is a linear combination with integer coefficients of products
\[
\prod_{h=1}^m \varphi(z_{\alpha(k_{1,h})}^{\epsilon(h)} \dots z_{\alpha(k_{r(h),h})}^{\epsilon(h)})
\]
where $m \ge 2$, the sets $\{k_{1,h} < \dots < k_{r(h),h}\}$ with $h = 1,\dots,m$ form a partition of $\{1,\dots,n\}$ and $\epsilon: \{1,\dots,m\} \to \{',''\}$.}

\bigskip
\noindent
{\bf {\em Proof.}} Expanding, we have
\[
\varphi((z'_{\alpha(1)} + z''_{\alpha(1)}) \dots (z'_{\alpha(n)} + z''_{\alpha(n)})) = \sum_{\epsilon \in \{',''\}^n} \varphi(z_{\alpha(1)}^{\epsilon(1)} \dots z_{\alpha(n)}^{\epsilon(n)}).
\]
The terms covered by ${\tilde Q}$ are those for which $\epsilon$ is not constant. Applying to each $\varphi(z_{\alpha(1)}^{\epsilon(1)} \dots z_{\alpha(n)}^{\epsilon(n)})$ Lemma~$5.2$ and Remark~$5.5$ the proposition follows. The fact that $m \ge 2$ is a consequence of the $\epsilon: \{1,\dots,n\} \to \{',''\}$ not being constant.\qed

\bigskip
\noindent
{\bf 5.7.} Having proved Proposition~$5.6$ it is now a standard application of Lie theory, like in \cite{12}, to derive an existence theorem for cumulants.

\bigskip
\noindent
{\bf Theorem.} {\em For each map $\alpha: \{1,\dots,n\} \to I \coprod J$ and $K \subset \{1,\dots,n\}$ let $\alpha K = \alpha \circ b_K$ where $b_K: \{1,\dots,|K|\} \to K$ is the increasing bijection. Then there is a universal polynomial
\[
\begin{aligned}
&R_{\alpha}(X_{\alpha K} \mid \emptyset \ne K \subset \{1,\dots,n\}) \\
&\quad = X_{\alpha\{1,\dots,n\}} + {\tilde R}_{\alpha}(X_{\alpha K} \mid \emptyset \ne K \subsetneqq \{1,\dots,n\})
\end{aligned}
\]
which is homogeneous of degree $n$ when $X_{\alpha K}$ is given degree $|K|$ and which has the cumulant property: if $z' = ((z'_i)_{i \in I},(z'_j)_{j \in J})$, $z'' = ((z''_i)_{i \in I},(z''_j)_{j \in J})$ is a bi-free pair of two-faced families of non-commutative random variables in $({\mathcal A},\varphi)$ and $z'+z'' = ((z'_i+z''_i)_{i \in I},(z'_j+z''_j)_{j \in J})$ then if $K = \{k_1 < \dots < k_r\} \subset \{1,\dots,n\}$ and $M'_{\alpha K} = \varphi(z'_{\alpha(k_1)} \dots z'_{\alpha(k_r)})$, $M''_{\alpha K} = \varphi(z''_{\alpha(k_1)} \dots z''_{\alpha(k_r)})$, $M'''_{\alpha K} = \varphi((z'_{\alpha(k_1)} + z''_{\alpha(k_1)}) \dots (z'_{\alpha(k_r)} + z''_{\alpha(k_r)}))$ we have
\[
\begin{aligned}
R_{\alpha}(M'_{\alpha K} \mid \emptyset \ne K \subset \{1,\dots,n\}) &+ R_{\alpha}(M''_{\alpha K} \mid \emptyset \ne K \subset \{1,\dots,n\}) \\
&= R_{\alpha}(M'''_{\alpha K} \mid \emptyset \ne K \subset \{1,\dots,n\}).
\end{aligned}
\]

The polynomial $R_{\alpha}$ of the above form satisfying the homogeneity requirement and the cumulant additivity property is unique.}

\bigskip
\noindent
{\bf {\em Proof.}} Let $\Pi'_n$ be set of non-empty subsets of $\{1,\dots,n\}$. We use Lemma~$5.6$ to turn ${\mathbb C}^{\alpha\Pi'_n}$ into a commutative algebraic group. Indeed, applying the lemma to $\alpha$ and to $\alpha|_K$ where $K \in \Pi'_n$, we have that there are polynomials
\[
\begin{aligned}
&Q_{\alpha K}(X_{\alpha L},Y_{\alpha L} \mid \alpha L \in \alpha\Pi'_n,L \subset K) \\
&\quad = X_{\alpha K} + Y_{\alpha K} + {\tilde Q}_{\alpha K}(X_{\alpha L},Y_{\alpha L})L \in \Pi'_n, L \subsetneqq K)
\end{aligned}
\]
so that if for $K = \{k_1 < \dots < k_r\} \subset \{1,\dots,n\}$, $M'_{\alpha K},M''_{\alpha K},M'''_{\alpha K}$ denote the moments of $z',z'',z'+z''$ corresponding to $\alpha|_K$ we have
\[
M'''_{\alpha K} = Q_{\alpha K}(M'_{\alpha L},M''_{\alpha L} \mid \alpha L \in \alpha\Pi'_n,\ L \subset K\}.
\]
Let $\boxplus\boxplus_n$ denote the operation on ${\mathbb C}^{\alpha \Pi'_n}$ defined by
\[
(\xi_{\alpha K} \mid \alpha K \in \alpha\Pi'_n) \boxplus\boxplus_n (\eta_{\alpha K} \mid \alpha K \in \alpha\Pi'_n) = (\zeta_{\alpha K} \mid \alpha K \in \alpha\Pi'_n)
\]
where $\zeta_{\alpha K} = Q_{\alpha K}(\xi_{\alpha L},\eta_{\alpha L} \mid \alpha L \in \alpha\Pi'_n,\ L \subset K)$. Since any system of complex numbers $(\xi_{\alpha K} \mid \alpha K \in \alpha\Pi'_n)$ are the corresponding moments for some two-faced family, we infer that $\boxplus\boxplus_n$ is commutative, associate and that $0$ is a neutral element. On the other hand the ``triangular'' feature of the algebraic law defining $\boxplus\boxplus_n$ shows the existence of inverses for the operation (the coordinates $\eta_{\alpha K}$ of the inverse of $(\xi_{\alpha K})_{\alpha K in \alpha\Pi'_n}$ if determined for $|K| < $ are then given for $|K| = m$ by
\[
\eta_{\alpha K} = -\xi_{\alpha K} - {\tilde Q}_{\alpha K}(\xi_{\alpha L},\eta_{\alpha L} \mid \alpha L \in \alpha\Pi'_n,\ L \subsetneqq K)
\]
with the observation that $|L| < m$). The exponential map $\exp: {\mathbb C}^{\alpha\Pi'_n} \to {\mathbb C}^{\alpha\Pi'_n}$ (the first ${\mathbb C}^{\alpha\Pi'_n}$ is the Lie algebra and the second is the group) one easily gets that
\[
\exp(\xi_{\alpha K} \mid \alpha K \in \alpha\Pi'_n) = (E_{\alpha K}(\xi_{\alpha L} \mid \alpha L \in \alpha\Pi'_n,\ \alpha L \subset \alpha K) \mid \alpha K \in \alpha \Pi'_n)
\]
where if $H \subset K$
\[
\frac {\partial E_{\alpha K}}{\partial \xi_{\alpha H}} = \frac {\partial Q_{\alpha K}}{\partial \xi_{\alpha H}} (0,E_{\alpha L} \mid \alpha L \in \alpha\Pi'_n,\ L \subset K)
\]
and $E_{\alpha K}(0) = 0$ Since $\frac {\partial Q_{\alpha K}}{\partial \xi_{\alpha K}} = 1$ one easily gets recurrently that $E_{\alpha K}$ is a polynomial in the $\xi_{\alpha L}, \emptyset \ne L \subset K$ and
\[
E_{\alpha K}(\xi_{\alpha L} \mid \alpha L \in \alpha\Pi'_n,\ L \subset K) = \xi_{\alpha K} + {\tilde E}_{\alpha K}(\xi_{\alpha L} \mid \alpha L \in \alpha\Pi'_n,\ L \subsetneqq K).
\]
Moreover one checks recurrently that $E_{\alpha K}$ is homogeneous of degree $|K|$ when $\xi_{\alpha L}$ is assigned degree $|L|$. Since the exponential $E$ has the homomorphism property
\[
\begin{aligned}
&(E_{\alpha K}(\xi_{\alpha L} \mid \alpha L \in \alpha\Pi'_n,\ L \subset K)) \boxplus\boxplus_n (E_{\alpha K}(\eta_{\alpha L} \mid \alpha L \in \alpha\Pi'_n,\ L \subset K)) \\
&\quad = E_{\alpha K}(\xi_{\alpha L} + \eta_{\alpha L} \mid \alpha L \in \alpha\Pi'_n,\ L \subset K)
\end{aligned}
\]
and $E_{\alpha K} = \xi_{\alpha K} + {\tilde E}_{\alpha K}$ implies that $E$ is a bijection, there is an inverse $(R_{\alpha K})_{\alpha K \in \alpha\Pi'_n}$ which is an homomorphism of $({\mathbb C}^{\alpha\Pi'_n},\boxplus\boxplus_n)$ and $({\mathbb C}^{\alpha\Pi'_n}),+)$ and has the properties that $R_{\alpha K}$ is homogeneous of degree $|K|$, depends only on variables $\xi_{\alpha L}$ with $L \subset K$ and $\xi_{\alpha K}$ has coefficient $1$ in $R_{\alpha K}$. Clearly taking $R_{\alpha} = R_{\alpha\{1,\dots,n\}}$ we get the polynomial with the desired properties.

Finally, the uniqueness of $R_{\alpha}$ can be seen as follows.

Since $R_{\alpha}$ is a differentiable homomorphism of $({\mathbb C}^{\alpha\Pi'_n},\boxplus\boxplus_n)$ into $({\mathbb C},+)$ it is completely determined by its differential at the neutral element, which is $(\xi_{\alpha K})_{\alpha K \in \alpha\Pi'_n} \to \xi_{\alpha\{1,\dots,n\}}$ since ${\tilde R}_{\alpha}$ has differential equal to $0$ at the origin, because of the homogeneity requirement, which implies there are no linear terms in the $\xi_{\alpha K}$'s in ${\tilde R}_{\alpha}$.\qed

\bigskip
\noindent
{\bf 5.8. Definition and Remark.} The polynomials $R_{\alpha}$ will be called bi-free cumulants. In case of a two-faced family $z = ((z_i)_{i \in I},(z_j)_{j \in J})$ in $({\mathcal A},\varphi)$ the {\em bi-free cumulants of $z$} are the numbers $R_{\alpha}(\varphi(z_{\alpha(k_1)} \dots z_{\alpha(k_r)}) \mid \{k_1 < \dots < k_r\} \subset \{1,\dots,n\})$. It will be often convenient to use various notations for these cumulants. If $z,\alpha$ are given and $\alpha(k) = a_k$ and $\mu_z$ is the distribution of $z$, we shall write for the corresponding bi-free cumulant $R_{a_1,\dots,a_n}(\mu_2)$, $R_{a_1\dots a_n}(z)$, $R_{\alpha}(\mu_z)$, $R_{\alpha}(z)$ or simply $R_{\alpha}$ or $R_{a_1\dots a_n}$ when it is clear to what two-faced distribution we refer.

\bigskip
Let us also record a further basic property of cumulants.

\bigskip
\noindent
{\bf 5.9. Proposition.} {\em For each map $\alpha: \{1,\dots,n\} \to I \coprod J$ there is a universal polynomial
\[
\begin{aligned}
M_{\alpha}(X_{\alpha K} \mid \emptyset \ne K \subset \{1,\dots,n\}) &= X_{\alpha\{1,\dots,n\}} \\
&+ {\tilde M}_{\alpha}(X_{\alpha K} \mid \emptyset \ne K \subsetneqq \{1,\dots,n\})
\end{aligned}
\]
which is homogeneous of degree $n$ when $X_{\alpha K}$ is assigned degree $(K)$ and such that
\[
M_{\alpha}(R_{\alpha K}^{(z)} \mid \emptyset \ne K \subset \{1,\dots,n\}) = \varphi(z_{\alpha(1)}\dots z_{\alpha(n)})
\]
if $z = ((z_i)_{i \in I},(z_j)_{j \in J})$ is a two-faced family in $({\mathcal A},\varphi)$.}

\bigskip
\noindent
{\bf {\em Proof.}} This is simply that the passage from moments has the ``triangular'' form which guarantees that it has an inverse of the same kind.\qed

\bigskip
\noindent
{\bf 5.10. Remark.} {\em If $\alpha: \{1,\dots,n\} \to I \coprod J$ is such that $\alpha(\{1,\dots,n\}) \subset I$ or $\alpha(\{1,\dots,n\}) \subset J$ then the bi-free cumulant $R_{\alpha}(z)$ reduces to a free cumulant of $(z_i)_{i \in I}$ or respectively of $(z_j)_{j \in J}$.} Indeed, for left or right variables bi-freeness and freeness are the same requirement and therefore in this case the uniqueness conditions for the free and bi-free cumulant are the same (Theorem~$5.7$ and \cite{9}, \cite{10}, \cite{12}).

\section{Examples}
\label{sec6}

\noindent
{\bf 6.1.} {\em Left and right regular representations of a free product of groups combined.} Let $(G)_{i \in I}$ be a family of discrete groups, let $G = {\underset{i \in I}{*}} G_i$ be their free product and let $({\mathbb C}[G_i],\tau_i)$, $({\mathbb C}[G],\tau)$ be the group algebras with their von~Neumann traces. Let further $\varphi: {\mathcal L}({\mathbb C}[G]) \to {\mathbb C}$ be the functional $\varphi(T) = \tau(T\delta_e)$, where $(\delta_g)_{g \in G}$ is the canonical basis of ${\mathbb C}[G]$ and let $L {\mathbb C}[G] \to {\mathcal L}({\mathbb C}[G])$ and $R: {\mathbb C}[G]^{op} \to {\mathcal L}({\mathbb C}[G])$ be the homomorphisms corresponding to the left and right regular representations, that is $L(\delta_h)\delta_g = \delta_{hg}$, $R(\delta_h)\delta_g = \delta_{gh}$. Viewing $G_i$ as a group of $G$, let ${\tilde L}_i = L \mid {\mathbb C}[G_i]$, ${\tilde R}_i = R \mid {\mathbb C}[G_i]^{op}$. {\em The family of pairs of faces $(({\mathbb C}[G_i],{\tilde L}_i)$, $({\mathbb C}[G_i]^{op},{\tilde R}_i))_{i \in I}$ is bi-free in $({\mathcal L}({\mathbb C}[G]),\varphi)$.} To see this, recall the identification of $({\mathbb C}[G],\ker \tau,\delta_e)$ with ${\underset{i \in I}{*}} ({\mathbb C}[G_i],\ker \tau,\delta_e)$ which sends $\delta_g$, $g \in G\backslash\{e\}$ to $\delta_{g_1} \otimes \dots \otimes \delta_{g_n} \in \ker \tau_{i_1} \otimes \dots \otimes \ker \tau_{i_n}$ when $g = g_1 \dots g_n$, $g_k \in G_{i_k}\backslash\{e\}$, $i_1 \ne i_2 \ne \dots \ne i_n$ is the reduced word for $g$. The maps $V_i$ and $W_i$ correspond essentially to the obvious identifications
\[
{\mathbb C}[G] \simeq {\mathbb C}[G_i] \otimes {\mathbb C}[G_i\backslash G] \simeq {\mathbb C}[G/G_i] \otimes {\mathbb C}[G_i].
\]
If $L_i,R_i$ are the left and right regular representation homomorphisms of ${\mathbb C}[G_i]$ and ${\mathbb C}[G_i]^{op}$ into ${\mathcal L}({\mathbb C}[G_i])$, then clearly if $h \in G_i$ we have $\lambda_i(L_i(\delta_h)) = {\tilde L}_i(\delta_h) = L(\delta_h)$ and $\rho_i(R_i(\delta_h)) = {\tilde R}_i(\delta_h) = R(\delta_h)$ and the expectation functional $\varphi$ is precisely the functional corresponding to the free product of vector spaces with specified state vector structure of $({\mathbb C}[G],\ker \tau,\delta_e)$, which proves the bi-freeness.

\bigskip
This example has clearly a Hilbertian version, which is generalized by the next example.

\bigskip
\noindent
{\bf 6.2.} {\em The standard form of a free product of von~Neumann algebras with faithful normal trace-state.} Let $(M_i,\tau_i)$, $i \in I$ be von~Neumann algebras with faithful normal trace state and consider their $W^*$-free product $(M,\tau) = {\underset{i \in I}{*}}(M_i,\tau_i)$. On $B(L^2(M,\tau))$ we denote by $\varphi$ the state $\varphi(T) = \langle T 1,1\rangle$ and by $L: M \to B(L^2(M,\tau))$, $R: M^{op} \to B(L^2(M,\tau))$ the left and right regular representations $L(m)h = mh$, $R(m)h = hm$, $h \in L^2(M,\tau)$. Viewing $M_i,M_i^{op}$ as subalgebras of $M,M^{op}$ let ${\tilde L}_i = L|_{M_i}: M_i \to B(L^2(M,\tau))$, ${\tilde R}_i = R|_{M_i^{op}}: M_i^{op} \to B(L^2(M,\tau))$ and consider also $L_i: M_i \to B(L^2(M_i,\tau_i))$, $R_i: M_i^{op} \to B(L^2(M_i,\tau_i))$ the left and right regular representations which are the restrictions of ${\tilde L}_i$ and ${\tilde R}_i$ to $L^2(M_i,\tau_i)$. {\em The family of pairs of faces $((M_i,{\tilde L}_i),(M_i^{op},{\tilde R}_i))_{i \in I}$ in $(B(L^2(M,\tau)),\varphi)$ is bi-free.} This is easily seen noticing that $(L^2(M,\tau),1) = {\underset{i \in I}{*}} (L^2(M_i,\tau_i),1)$ and $\lambda_i(L_i(m)) = {\tilde L}_i(m)$, $\rho_i(R_i(m)) = {\tilde R}_i(m)$ if $m \in M_i$.

\bigskip
\noindent
{\bf 6.3.} {\em Left and right creation and annihilation operators on the full Fock space combined.} Let ${\mathcal H}$ be a complex Hilbert space with orthonormal basis $(e_i)_{i \in I}$ and let ${\mathcal T}({\mathcal H}) = {\mathbb C} 1 \oplus \underset{n \ge 1}{\bigoplus} {\mathcal H}^{\otimes n}$ be the full Fock space on which left and right creation operators $l_i\zeta = e_i \otimes \zeta$, $r_i\zeta = \zeta \otimes e_i$ if $\zeta \in \underset{n \ge 1}{\bigoplus} {\mathcal H}^{\otimes n}$ and $l_i1 = r_i1 = e_i$. Their adjoints are the annihilation operators $l_i^*,r_i^*$. {\em Then the family of two-faced families $((l_i,l_i^*),(r_i,r_i^*))_{i \in I}$ is bi-free on $({\mathcal B}({\mathcal T}({\mathcal H}))$, $\omega_{{\mathcal H}})$ where $\omega_{{\mathcal H}}(\cdot) = \langle \cdot 1,1\rangle$.} Indeed, this follows from the natural identification $({\mathcal T}({\mathcal H}),1) = {\underset{i \in I}{*}}({\mathcal T}({\mathbb C} e_i),1)$ under which $\lambda_i(l_i \mid {\mathcal T}({\mathbb C} e_i)) = l_c$, $\lambda_i(l_i^* \mid {\mathcal T}({\mathbb C} e_i)) = l_i^*$, $\rho_i(r_i \mid {\mathcal T}({\mathbb C} e_i)) = r_i$, $\rho_i(r_i^* \mid {\mathcal T}({\mathbb C} e_i)) = r^*$. Using $2.12$ and Proposition~$2.13$ we can give a generalization of the previous bi-freeness fact. Assume ${\mathcal H} = \underset{k \in K}{\bigoplus} {\mathcal H}_k$ and consider creation and annihilation operators $l(h)\zeta = h \otimes \zeta$, $r(h)\zeta = \zeta \otimes h$, $l(h)1 = r(h)1 = h$ and $l^*(h) = l(h)^*$, $r^*(h) = r(h)^*$. {\em Then the family of pairs of faces $(C^*(l({\mathcal H}_k)),C^*(r({\mathcal H}_k)))_{k \in K}$ is bi-free in $({\mathcal B}({\mathcal T}({\mathcal H})),\omega_{{\mathcal H}})$.} Indeed, we can find an orthonormal basis $(e_i)_{i \in I}$ in ${\mathcal H}$ for which there is a partition $I = \underset{k \in K}{\coprod} I_k$ so that $(e_i)_{i \in I_k}$ is an orthonormal basis of ${\mathcal H}_k$. Using $2.12$ and Proposition~$2.13$ together with preservation of bi-freeness under norm-closure we get the desired bi-freeness fact.

\newpage
\section{The Algebraic Bi-free Central Limit Distributions and Theorem}
\label{sec7}

\noindent
{\bf 7.1.} At this stage, like in our first free probability paper \cite{12}, having proved an existence theorem for cumulants and having enough examples of distributions which can be central limit distributions, we can prove an algebraic bi-free central limit theorem for the convergence of moments for two-faced families of non-commutative random variables. The first step will be to identify the bi-free cumulants of degree $1$ and $2$ and to present the candidates for central limit distributions.

\bigskip
\noindent
{\bf 7.2. Lemma.} {\em If $z = ((z_i)_{i \in I},(z_j)_{j \in J})$ is a two-faced family of non-commutative random variables in $({\mathcal A},\varphi)$ and $a,b \in I \coprod J$ the bi-free cumulants of degree $\le 2$ are given by the formulae
\[
\begin{aligned}
R_a(\mu_z) &= \varphi(z_a), \\
R_{ab}(\mu_z) &= \varphi(z_az_b) - \varphi(z_a)\varphi(z_b).
\end{aligned}
\]
}

\bigskip
\noindent
{\bf {\em Proof.}} In view of the uniqueness of the cumulant polynomials defined in Theorem~$5.7$ we see that the polynomials in the statement of the Lemma have the appropriate form and homogeneity properties, so the proof reduces to checking that they satisfy the addition property of cumulants. Thus consider $z^{\epsilon} = ((z_i^{\epsilon})_{i \in I}, (z_j^{\epsilon})_{j \in J})$, $\epsilon \in \{',''\}$ a bi-free in $({\mathcal A},\varphi)$. For degree~$1$, $\varphi((z'+z'')_a) = \varphi(z'_a) + \varphi(z''_a)$ obviously. For the degree~$2$ formulae, remark first that if both $a,b \in I$ or if both $a,b \in J$, then bi-freeness reduces to freeness of this pair of indices and the bi-free cumulants $R_{ab}$ must be free cumulants and these are known formulae. Thus, we are left with the mixed cases $a \in I$, $b \in J$ and $a \in J$, $b \in I$. To simplify matters remark that $\varphi((z_a+c_a1)(z_v+c_b1)) - \varphi(z_a+c_a1)\varphi(z_b+c_b1) = \varphi(z_az_b) - \varphi(z_a)\varphi(z_b)$ and that it suffices to check the addition property in case $\varphi(z'_a) = \varphi(z'_b) = \varphi(z''_a) = \varphi(z''_b) = 0$. We must then show that
\[
\begin{aligned}
0 &= \varphi((z'_a+z''_a)(z'_b+z''_b)) - \varphi(z'_az'_b) - \varphi(z''_az''_b) \\
&= \varphi(z'_az''_b) + \varphi(z''_az'_b).
\end{aligned}
\]
Thus it suffices to show that $\varphi(z'_az''_b) = \varphi(z''_az'_b) = 0$ under these assumptions. In view of the definition of bi-freeness this reduces to statements about certain left and right operators on free products of vector spaces with specified state vectors. Thus, let $({\mathcal X}^{\epsilon},{\mathcal X}^{\epsilon}_0,\xi^{\epsilon})$, $\epsilon \in \{',''\}$ be vector spaces with specified state vectors and $({\mathcal X},{\mathcal X}_0,\xi) = ({\mathcal X}',{\mathcal X}'_0,\xi') * ({\mathcal X}'',{\mathcal X}''_0,\xi'')$. Let further $T_a^{\epsilon},T_b^{\epsilon} \in {\mathcal L}({\mathcal X}^{\epsilon})$, $\varphi_a^{\epsilon} = \varphi_{\epsilon}(T_b^{\epsilon}) = 0$, that is $T_a^{\epsilon}\xi^{\epsilon} \in {\mathcal X}_0^{\epsilon}$, $T_b^{\epsilon}\xi^{\epsilon} \in {\mathcal X}_0^{\epsilon}$. We must show that
\[
\begin{aligned}
\varphi(\lambda,(T'_a)\rho,,(T''_b)) &= \varphi(\rho,,(T''_b)\lambda,(T'_a)) \\
&= \varphi(\lambda,,(T''_a)\rho,(T'_b) \\
&= \varphi(\rho,(T'_b)\lambda,,(T''_a)) = 0
\end{aligned}
\]
which is clearly so since
\[
\begin{aligned}
\lambda,(T'_a)\rho,,(T''_b)\xi &= \rho,,(T''_a)\lambda,(T'_a)\xi \\
&= T'_a\xi' \otimes T''_b\xi'' \in {\mathcal X}'_0 \otimes {\mathcal X}''_0 \\
\lambda,,(T''_a)\rho,(T'_b)\xi &= \rho,(T'_b)\lambda,,(T''_a)\xi \\
&= T''_a\xi'' \otimes T'_b\xi' \in {\mathcal X}''_0 \otimes {\mathcal X}'_0.
\end{aligned}
\]
\qed

\bigskip
\noindent
{\bf 7.3. Definition.} A two-faced family $z = ((z_i)_{i \in I},(z_j)_{j \in J})$ in $({\mathcal A},\varphi)$ has a {\em bi-free central limit distribution} (or {\em centered bi-free Gaussian distribution}) if its cumulants satisfy $R_{a_1\dots a_k}(\mu_z) = 0$ if $k \ge 3$ or $k = 1$.

\bigskip
\noindent
{\bf 7.4. Theorem.} {\em There is exactly one bi-free central limit distribution $\gamma_C: {\mathbb C}\langle Z_k \mid k \in I \coprod J\rangle \to {\mathbb C}$ for a given matrix $C = (C_{kl})_{k,l \in I \coprod J}$ with complex entries so that
\[
\gamma_C(Z_kZ_l) = C_{kl},\ k,l \in I \coprod J.
\]
If $h,h^*: I \coprod J \to {\mathcal H}$ are maps into the Hilbert space ${\mathcal H}$ and with the notation of $6.3$ we define
\[
\begin{aligned}
z_i &= l(h(i)) + l^*(h^*(i)) \text{ if $i \in I$,} \\
z_j &= r(h(j)) + r^*(h^*(j)) \text{ if $j \in J$,}
\end{aligned}
\]
then $z = ((z_i)_{i \in I},(z_j)_{j \in J})$ has a bi-free central limit distribution $\gamma_C$ where $C_{kl} = \langle h(l),h^*(l)\rangle$. Every bi-free central limit distribution in case $I,J$ are finite, can be obtained in this way.}

\bigskip
\noindent
{\bf {\em Proof.}} Since a distribution is determined by its cumulants, a central limit distribution is completely determined by its second order moments in view of the formulae in Lemma~$7.2$. Also because any set of cumulants is the set of cumulants of the algebraic distribution of some bi-free family we get all central limit distributions in this way.

For the remaining part of the theorem remark that $z_l1 = h(l)$ and $z^*_k1 = h^*(k)$ so that
\[
\langle z_kz_l1,1\rangle = \langle z_l1,z_k^*1\rangle = \langle h(l),h^*(k)\rangle.
\]
Also clearly $\langle z_k1,1\rangle = 0$. Hence to show that $z$ has distribution $\gamma_C$ it will suffice to show that $R_{a_1\dots a_k}(\mu_z) = 0$ when $k \ge 3$. This follows from the usual considerations based on homogeneity. In ${\mathcal H} \oplus {\mathcal H}$ let $h'(k) = h(k) \oplus 0$, $h''(k) = 0 \oplus h(k)$, $h^*{}'(k) = h^*(k) \oplus 0$, $h^*{}''(k) = 0 \oplus h^*(k)$ and define $z',z''$ using the $h',h^*{}'$ and respectively $h'',h^*{}''$ instead of $h,h^*$ vectors. Then $z',z''$ are bi-free and $2^{-1/2}(z'+z'')$ has the same distribution like $z$ since the map ${\mathcal H} \ni h \to 2^{-1/2}(h \oplus h) \in {\mathcal H} \oplus {\mathcal H}$ is isometric. This then gives $R_{a_1\dots a_k}(z) = R_{a_1\dots a_k}(2^{-1/2}(z'+z'')) = 2^{-k/2}(R_{a_1\dots a_k}(z')) + 2^{-k/2}R_{a_1\dots a_k}(z'') = 2^{1-k/2}R_{a_1\dots a_k}(z)$. Since $2^{1-k/2} \ne 1$ if $k \ne 2$, we get $R_{a_1\dots a_k}(z) = 0$ if $k \ne 2$.\qed

\bigskip
\noindent
{\bf 7.5. Definition.} A $*$-distribution
\[
\varphi: {\mathbb C}\langle Z_k,Z^*_k \mid k \in I \coprod J\rangle \to {\mathbb C}
\]
is a {\em bi-free central limit $*$-distribution}, if it is a bi-free central limit distribution for $((Z_i,Z_i^*)_{i \in I},(Z_j,Z^*_j)_{j \in J})$ and satisfies the positivity condition $\varphi(P^*P) \ge 0$ for $P \in {\mathbb C}\langle Z_k,Z^*_k \mid k \in I \coprod J\rangle$.

\bigskip
\noindent
{\bf 7.6. Theorem.} {\em A bi-free central limit distribution
\[
\varphi: {\mathbb C}\langle Z_k,Z^*_k \mid k \in I \coprod J\rangle \to {\mathbb C}
\]
is a bi-free central limit $*$-distribution iff the matrix of second order moments
\[
C = \left( \begin{matrix}
(\varphi(Z^*_kZ_l))_{k,l \in K} & (\varphi(Z_kZ_l))_{k,l \in K} \\
(\varphi(Z^*_kZ^*_l))_{k,l \in K} & (\varphi(Z_kZ^*_l))_{k,l \in K}
\end{matrix} \right)
\]
where $K = I \coprod J$ is $\ge 0$.

If $K = I \coprod J$ is finite and $h(k),h^*(l) \in {\mathcal H}$ are vectors in the Hilbert space ${\mathcal H}$ such that
\[
C = \left( \begin{matrix}
(\langle h(l),h(k)\rangle)_{k,l \in K} & (\langle h(l),h^*(k)\rangle )_{k,l \in K} \\
(\langle h^*(l),h(k)\rangle)_{k,l \in K} & (\langle h^*(l),h^*(k)\rangle )_{k,l \in K}
\end{matrix} \right)
\]
then the $*$-distribution with respect to the vacuum vector $1 \in T({\mathcal H})$ of the two-faced family $z = ((z_i)_{i \in I},(z_j)_{j \in J})$, where $z_k = l(h(i)) + l^*(h^*(i))$, $i \in I$, $z_j = r(h(j)) + r^*(h^*(j))$, $j \in J$ is precisely the bi-free central limit $*$-distribution with second order moments matrix $C$.}

\bigskip
\noindent
{\bf {\em Proof.}} It is easily observed that it suffices to prove the assertions under the assumption that $K = I \coprod J$ is finite. Clearly if $\varphi$ is positive $C \ge 0$. In view of the uniqueness of a bi-free central limit distribution with given second order moments (Theorem~$7.5$), that $C \ge 0$ implies $\varphi \ge 0$ will follow from the second assertion about $\varphi$ being the $*$-distribution of Hilbert space operators. Remark that for every positive matrix $C$ there are vectors $h(k),h^*(k) \in {\mathcal H}$ so that it is a matrix of scalar products as above. Further, the two-faced family $(z,z^*) = ((z_i,z^*_i)_{i \in I},(z_j,z^*_j)_{j \in J})$ has a bi-free central limit distribution by Theorem~$7.4$. Indeed, if $h(k),h^*(k),k \in K$ yield the given matrix $C$ then $(z,z^*)$ with
\[
\begin{aligned}
z_i &= l(h^*(i)) + l^*(h^*(i)), \\
z^*_i &= l(h^*(i)) + l^*(h(i)), \\
z_j &= r(h(j)) + r^*(h^*(j)), \\
z^*_j &= r(h^*(j)) + r^*(h(j))
\end{aligned}
\]
has bi-free central limit distribution and matrix $C$ of second order moments.\qed

\bigskip
\noindent
{\bf 7.7. Definition.} A distribution
\[
\varphi: {\mathbb C}\langle X_k \mid k \in I \coprod J\rangle \to {\mathbb C}
\]
is a hermitian bi-free central limit distribution (w.r.t.\ $I$ left and $J$ right index sets) and $\varphi \ge 0$ when ${\mathbb C}\langle X_k \mid k \in I \coprod J\rangle$ is endowed with the $*$-algebra structure in which $X_k = X^*_k$ for all $k \in I \coprod J$.

\bigskip
\noindent
{\bf 7.8. Theorem.} {\em A bi-free central limit distribution $\varphi: {\mathbb C}\langle X_k \mid k \in I \coprod J\rangle \to {\mathbb C}$ is hermitian iff the matrix of second order moments $C = (\varphi(X_kX_l))_{k,l \in I \coprod J}$ is $\ge 0$.}

\bigskip
If $K = I \coprod J$ is finite and $h(k) \in {\mathcal H}$ are such that
\[
C = (\langle h(l),h(k)\rangle)_{h,l \in K}
\]
then the two-faced hermitian family $x = ((x_i)_{i \in I},(x_j)_{j \in J})$ on $({\mathcal T}({\mathcal H}),1)$ where $x_i = l(h(i)) + l^*(h(i))$, $i \in I$ and $x_j = r(h(j)) + r^*(h(j))$ has precisely the hermitian bi-free central limit distribution with second order moments matrix $C$ as distribution.

\bigskip
\noindent
{\bf {\em Proof.}} Like in the proof of the preceding theorem it suffices to carry out the proof when $K$ is finite and clearly $\varphi \ge 0$ implies $C \ge 0$. Again the converse $C \ge 0 \Rightarrow \varphi \ge 0$ will follow from the uniqueness of $\gamma_C$ in Theorem~$7.4$ and the second assertion of the theorem which produces a hermitian two-faced family with given matrix $C$. The last assertion follows easily from Theorem~$7.4$ applied to vectors $h(k) = h^*(k)$, $k \in K$.\qed

\bigskip
\noindent
{\bf 7.9. Theorem (Bi-free Algebraic Central Limit Theorem).} {\em Let $z^{(n)} = ((z_i^{(n)})_{i \in I},(z_j^{(n)})_{j \in J})$, $n \in {\mathbb N}$ be a bi-free sequence of two-faced families in $({\mathcal A},\varphi)$, such that
\begin{itemize}
\item[{\rm (i)}] $\varphi(z_k^{(n)}) = 0$, $k \in I \coprod J$
\item[{\rm (ii)}] $\underset{n \in {\mathbb N}}{\sup} |\varphi(z_{k_1}^{(n)} \dots z_{k_n}^{(n)})| = D_{k_1\dots k_m} < \infty$ for every $k_1,\dots,k_m \in I \coprod J$
\item[{\rm (iii)}] $\underset{N \to \infty}{\lim} N^{-1} \underset{1 \le n \le N}{\sum} \varphi(z_k^{(n)}z_l^{(n)}) = C_{kl}$ for every $k,l \in I \coprod J$.
\end{itemize}
}

\bigskip
Then letting $S_N = ((S_{N,i})_{i \in I},(S_{N,j})_{j \in J})$, $S_{N,k} = N^{-1/2} \underset{1 \le n \le N}{\sum} z_k^{(n)}$, $k \in I \coprod J$ and letting $\gamma_C$ be the bi-free central limit distribution in Theorem~$7.4$ with $C = (C_{k,l})_{k,l \in I,J}$, we have
\[
\lim_{N \to \infty} \mu_{S_N}(P) = \gamma_C(P)
\]
for every $P \in {\mathbb C}\langle Z_k \mid k \in I \coprod J\rangle$.

\bigskip
\noindent
{\bf {\em Proof.}} Since moments are polynomials of cumulants by Proposition~$5.9$, also in the bi-free context it suffices to prove that
\[
\lim_{N \to \infty} R_{k_1\dots k_m}(\mu_{S_N}) = R_{k_1 \dots k_m}(\gamma_C).
\]
In view of the definition of $\gamma_C$, this means
\[
\begin{aligned}
&\lim_{N \to \infty} R_{kl}(\mu_{S_N}) = C_{kl} \text{ and} \\
&\lim_{N \to \infty} R_{k_1\dots k_m}(\mu_{S_N}) = 0 \text{ if $m \ne 2$.}
\end{aligned}
\]
In view of Lemma~$7.2$ and of the properties of cumulants, we have $R_k(\mu_{S_N}) = \varphi(S_N) = 0$ and
\[
\begin{aligned}
R_{k_1\dots k_m}(\mu_{S_N}) &= \sum_{1 \le n \le N} R_{k_1\dots k_m}(N^{-1/2}\mu_{z^{(n)}}) \\
&= N^{-m/2} \sum_{1 \le n \le N} R_{k_1\dots k_m}(\mu_{z^{(n)}}).
\end{aligned}
\]
If $m = 2$ we have
\[
\begin{aligned}
R_{kl}(\mu_{S_N}) &= N^{-1} \sum_{1 \le n \le N} R_{kl}(\mu_{z^{(n)}}) \\
&= N^{-1} \sum_{1 \le n \le N} \varphi(z_k^{(n)}z_l^{(n)})
\end{aligned}
\]
and hence in view of (iii)
\[
\lim_{N \to \infty} R_{kl}(\mu_{S_N}) = C_{kl}.
\]
Since the bi-free cumulants are universal polynomials in the moments we infer from (ii) that
\[
\sup_{n \in {\mathbb N}} |R_{k_1\dots k_m}(z^{(n)})| = B_{k_1\dots k_m} < \infty
\]
and hence if $m \ge 3$
\[
\lim_{N \to \infty} R_{k_1 \dots k_m}(\mu_{S_N}) = 0.
\]
\qed

\bigskip
\noindent
{\bf 7.10. Corollary.} {\em If in Theorem~$7.9$ the non-commutative probability space $({\mathcal A},\varphi)$ is a $C^*$-probability space and the $z^{(n)}$ are hermitian, that is $z_k^{(n)} = z_k^{(n)*}$ for all $k \in I \coprod J$, then the bi-free central limit distribution $\gamma_C$ in the theorem is a bi-free hermitian central limit distribution.}

\bigskip
\noindent
{\bf 7.11. Remark.} Modifying the assumptions of Theorem~$7.9$ to involve a $C^*$-probability space $({\mathcal A},\varphi)$ and the $*$-distributions of the $z^{(n)}$ instead of their distributions one infers from the theorem applied to the $(z^{(n)},z^{(n)*})$ a convergence in $*$-distribution of the $S_N$ to a bi-free central limit $*$-distribution. We leave the details as an exercise to the reader.

\section{Bi-freeness with Amalgamation Over $B$ and $B-B$ Non-commutative Probability Spaces}
\label{sec8}

\noindent
{\bf 8.1.} We will sketch how to generalize the definition of bi-freeness when ${\mathbb C}$ is replaced by a unital algebra $B$ over ${\mathbb C}$ and we will introduce the $B-B$ non-commutative probability spaces, which we will use for this instead of the usual non-commutative probability spaces over $B$. We will pursue elsewhere, the generalization along these lines of other parts of the present paper.

\bigskip
\noindent
{\bf 8.2.} If ${\mathcal X}$ is a vector space over ${\mathbb C}$, which is a $B-B$ bimodule, we will denote by ${\mathcal L}({\mathcal X})$, ${\mathcal L}_l({\mathcal X})$ and ${\mathcal L}_r({\mathcal X})$ the algebra of all operators on ${\mathcal X}$ and the two subalgebras of operators which are $B$-linear with respect to the left and respectively right $B$-module structures on $B$. The left and right $B$-module structures on ${\mathcal X}$ yield a homomorphism of $B \otimes B^{op}$ into ${\mathcal L}({\mathcal X})$ and ${\mathcal L}_l({\mathcal X}),{\mathcal L}_r({\mathcal X})$ are the commutants of the homomorphic images of $B \otimes 1$ and $1 \otimes B^{op}$. We will be interested in bi-modules with specified ``state-vector'', that is bimodules ${\mathcal X} = B \oplus \overset{\circ}{{\mathcal X}}$, where $\overset{\circ}{{\mathcal X}}$ is also a bimodule. If $T \in {\mathcal L}({\mathcal X})$ we define $p_{{\mathcal X}}(T)$ (abbreviated $p(T)$ sometimes), the element in $B$ so that $T(1 \oplus O) \in p_{{\mathcal X}}(T) \oplus \overset{\circ}{{\mathcal X}}$. The structure which arises on ${\mathcal L}({\mathcal X})$ is abstracted in the next definition.

\bigskip
\noindent
{\bf 8.3. Definition.} A $B-B$ non-commutative probability space is a triple $({\mathcal A},p,\epsilon)$ where ${\mathcal A}$ is a unital algebra over ${\mathbb C}$, $\epsilon: B \otimes B^0 \to {\mathcal A}$ is a unital homomorphism so that $\epsilon \mid B \otimes 1$ and $\epsilon \mid 1 \otimes B^{op}$ are injective and $p: {\mathcal A} \to B$ is a linear map so that
\[
p(\epsilon(b_1 \otimes 1)a\epsilon(1 \otimes b_2)) = b_1p(a)b_2
\]
and in particular $(p \circ \epsilon)(b_1 \otimes b_2) = b_1b_2$ if $a \in {\mathcal A}$, $b_1,b_2 \in B$. The commutants ${\mathcal A}_l$, ${\mathcal A}_r$ of $\epsilon(B \otimes 1)$ and $\epsilon(1 \otimes B^{op})$ in ${\mathcal A}$, will be called the algebras of left and right variables.

\bigskip
\noindent
{\bf 8.4.} If ${\mathcal X}_i = B \oplus \overset{\circ}{{\mathcal X}}_i$, $i \in I$ are bi-modules with specified unit vectors, their free product ${\mathcal X} = B \oplus \overset{\circ}{{\mathcal X}}$ also denoted $({\mathcal X},\overset{\circ}{{\mathcal X}}) = {\underset{i \in I}{*_B}} ({\mathcal X}_i,\overset{\circ}{{\mathcal X}}_i)$ is defined by
\[
\overset{\circ}{{\mathcal X}} = \underset{n \ge 1}{\oplus} \underset{\underset{i_i \ne \dots \ne i_n}{i_j \in I}}{\bigoplus} \overset{\circ}{{\mathcal X}}_{i_1} \otimes_B \dots \otimes_B \overset{\circ}{{\mathcal X}}_{i_n}.
\]
We consider the isomorphisms of bimodules
\[
V_i : {\mathcal X} \to {\mathcal X}_i \otimes_B (B \oplus \underset{n \ge 1}{\oplus} \underset{\underset{i_1 \ne \dots \ne i_n}{i_j \in I,i \ne i_1}}{\bigoplus} \overset{\circ}{{\mathcal X}}_{i_1} \otimes \dots \otimes \overset{\circ}{{\mathcal X}}_{i_n})
\]
and
\[
W_i: {\mathcal X} \to (B \oplus \underset{n \ge 1}{\oplus} \underset{\underset{i_1 \ne \dots \ne i_n}{i_j \in I,i_n \ne i}}{\bigoplus} \overset{\circ}{{\mathcal X}}_{i_1} \otimes \dots \otimes \overset{\circ}{{\mathcal X}}_{i_n}) \otimes_B {\mathcal X}_i
\]
and the unital homomorphisms
\[
\rho_i: {\mathcal L}_l({\mathcal X}_i) \to {\mathcal L}_l({\mathcal X}),\lambda_i: {\mathcal L}_r({\mathcal X}_i) \to {\mathcal L}_r({\mathcal X})
\]
defined by
\[
\begin{aligned}
\rho_i(T) &= W_i^{-1}(I \otimes T)W_i, \\
\lambda_i(T) &= V_i^{-1}(T \otimes I)V_i.
\end{aligned}
\]

\bigskip
\noindent
{\bf 8.5. Definition.} If $({\mathcal A},p,\epsilon)$ is a $B-B$ non-commutative probability space, a {\em family of pairs of $B$-faces}, is a family $(C_i,D_i)_{i \in I}$ of unital subalgebras of ${\mathcal A}$ so that $\epsilon(B \otimes 1) \subset C_i \subset {\mathcal A}_l$, $\epsilon(1 \otimes B) \subset D_i \subset {\mathcal A}_r$ if $i \in I$.

The family of $B$-faces $(C_i,D_i)_{i \in I}$ is {\em bi-free over $B$} if there are $B-B$ bimodules ${\mathcal X}_i = B \oplus \overset{\circ}{{\mathcal X}}_i$ and homomorphisms $\gamma_i: C_i \to {\mathcal L}_r({\mathcal X}_i)$, $\delta_i: D_i \to {\mathcal L}_l({\mathcal X}_i)$ so that $\gamma_i(\epsilon(b \otimes 1))$, $\delta_i(\epsilon(1 \otimes b))$ are the linear operators of left and right multiplication by $b$ on ${\mathcal X}_i$ and on $({\mathcal X},\overset{\circ}{{\mathcal X}}) = {\underset{i \in I}{*_B}} ({\mathcal X}_i,\overset{\circ}{{\mathcal X}}_i)$ the following condition is satisfied: 

\noindent
if $c_k \in C_{i(k)}$, $d_k \in D_{i(k)}$, $1 \le k \le n$

\noindent
then:

\[
\begin{aligned}
&p_{{\mathcal X}}(\lambda_{i(1)}(\gamma_{i(1)}(c_1))\rho_{i(1)}(\delta_{i(1)}(d_1))\dots \lambda_{i(n)}(\gamma_{i(n)}(c_n))\rho_{i(n)}(\delta_{i(n)}(d_n))) \\
&\quad = p(c_1d_1c_2d_2 \dots c_nd_n).
\end{aligned}
\]

\end{document}